\documentclass{amsart}
\usepackage{graphicx}
\usepackage{amssymb}
\usepackage[all]{xy}
\xyoption{poly}

\newcommand{\Dkf}{\operatorname{Diff}^{k}_{\Phi}}

\newcommand{\FDens}{{}^{\Phi}\Omega}

\newcommand\ff{\operatorname{ff}}

\numberwithin{equation}{section}

\newcommand\paperbody%
        {}

\newtheorem{lemma}{Lemma}
\newtheorem{prop}{Proposition}
\newtheorem{corollary}{Corollary}
\newtheorem{prob}{Problem}
\newtheorem{theorem}{Theorem}

\theoremstyle{remark}
\newtheorem{defin}{Definition}
\newtheorem{remark}{Remark}

\newcommand\FT{{}^{\Phi}T}
\newcommand\fT{{}^{\phi}T}
\newcommand\fS{{}^{\phi}S}
\newcommand\fpi{{}^{\phi}\pi}
\newcommand\ftpi{{}^{\phi}\tilde\pi}

\newcommand\cFTs{{}^{\Phi}\overline{T}\kern-1pt{}^*}

\newcommand\sus{\operatorname{sus}}

\newcommand\ie{i\@.e\@.\ }

\newcommand\WF{\operatorname{WF}}

\newcommand\Char{\operatorname{Char}}

\newcommand\Charb{\operatorname{Char}_{\pa}}

\newcommand\Ells{\operatorname{Ell}_{\sigma}}
\newcommand\Ellb{\operatorname{Ell}_{\pa}}

\newcommand\Dens{\Omega}

\newcommand\CC{\mathbb C}

\newcommand\NN{\mathbb N}

\newcommand\RR{\mathbb R}
\renewcommand\SS{\mathbb S}

\newcommand\ZZ{\mathbb Z}
\newcommand\CIc{{\mathcal{C}}^{\infty}_c}

\newcommand\CI{{\mathcal{C}}^{\infty}}
\newcommand\CmIc{{\mathcal{C}}^{-\infty}_c}
\newcommand\CmI{{\mathcal{C}}^{-\infty}}

\newcommand\Diag{\operatorname{Diag}}
\newcommand\Diff[1]{\operatorname{Diff}^{#1}}
\newcommand\DiffF[1]{\operatorname{Diff}_{\Phi}^{#1}}
\newcommand\PsiF[1]{\Psi_{\Phi}^{#1}}
\newcommand\Psisc[1]{\Psi_{\text{sc}}^{#1}}
\newcommand\PsipF[1]{\Psi_{p\Phi}^{#1}}
\newcommand\Psif[1]{\Psi_{\phi}^{#1}}
\newcommand\Psis[2]{\Psi_{\sus(#1)}^{#2}}
\newcommand\Psisf[2]{\Psi_{\sus(#1)-\phi}^{#2}}
\newcommand\XF[1]{X_{\Phi}^{#1}}
\newcommand\CF{C_{\Phi}}
\newcommand\WFF{\operatorname{WF}_{\Phi}}
\newcommand\WFsc{\operatorname{WF}_{\text{sc}}}
\newcommand\WFFs{\operatorname{WF}_{\Phi,\sigma}}
\newcommand\WFFb{\operatorname{WF}_{\Phi,\pa}}

\newcommand\FN{{}^{\Phi}N}
\newcommand\cscT{{}^{\text{sc}}\overline{T}}
\newcommand\FSN{{}^{\Phi}SN}

\newcommand\cFN{{}^{\Phi}\overline N}
\newcommand\cFNs{{}^{\Phi}\overline N\kern-1pt{}^*}
\newcommand\FS{{}^{\Phi}S}

\newcommand\NF{N_{\Phi}}

\newcommand\Hom{\operatorname{Hom}}
\newcommand\Id{\operatorname{Id}}
\newcommand\Lap{\varDelta}

\newcommand\ci{${\mathcal{C}}^\infty$}
\newcommand\clos{\operatorname{cl}}

\newcommand\dCI{\dot{\mathcal{C}}^{\infty}}

\newcommand\dCmI{\dot{\mathcal{C}}^{-\infty}}
\newcommand\ha{\frac{1}{2}}

\newcommand\pa{\partial}
\newcommand\paX{\partial X}
\newcommand\pbX{Y}

\newcommand\restrictedto{\upharpoonright}

\newcommand\supp{\operatorname{supp}}

\newcommand\Fa{\operatorname{Fa}}

\newcommand\longhookrightarrow{\hookrightarrow}
\newcommand\longtwoheadrightarrow{\twoheadrightarrow}

\newcommand\Mand{\text{ and }}

\newcommand\Mfor{\text{ for }}

\newcommand\Mforsome{\text{ for some }}

\newcommand\Mif{\text{ if }}

\newcommand\Min{\text{ in }}

\newcommand\Mst{\text{ s.t. }}

\newcommand\Mwith{\text{ with }}

\title[Manifolds with fibred boundaries]
{Pseudodifferential operators on manifolds with fibred boundaries}

\author[Rafe Mazzeo]{Rafe Mazzeo$^{\dag}$}
\address{Department of Mathematics, Stanford University}
\email{mazzeo@math.stanford.edu}
\author[Richard Melrose]{Richard B. Melrose$^{*}$}
\address{Department of Mathematics, MIT}
\email{rbm@math.mit.edu}

\thanks{This research was supported in part by the National Science
Foundation under an NYI and grants DMS-9626382 ($^{\dag}$) and 
DMS-9306389 ($^{\S}$). Manuscript available from
http:{\scriptsize/}{\scriptsize/}www-math.mit.edu{\scriptsize/}{$\sim$}rbm{\scriptsize/}rbm-home.html}

\dedicatory{Respectfully dedicate to Professor M. Sato on the occasion of
his $70$th birthday}

\begin{document}

\begin{abstract} Let $X$ be a compact manifold with boundary. Suppose that
the boundary is fibred, $\phi:\pa X\longrightarrow Y,$ and let $x\in\CI(X)$
be a boundary defining function. This data fixes the space of `fibred cusp'
vector fields, consisting of those vector fields $V$ on $X$ satisfying
$Vx=O(x^2)$ and which are tangent to the fibres of $\phi;$ it is a Lie
algebra and $\CI(X)$ module. This Lie algebra is quantized to the `small
calculus' of pseudodifferential operators $\PsiF*(X).$ Mapping
properties including boundedness, regularity, Fredholm condition and symbolic
maps are discussed for this calculus. The spectrum of the Laplacian of an
`exact fibred cusp' metric is analyzed as is the wavefront set associated
to the calculus.
\end{abstract}

\maketitle

\section*[Introduction]{Introduction\label{S.I}}

Algebras of pseudodifferential operators can be used to investigate local
regularity of solutions to partial differential equations and to relate
such local matters to more global properties. On a compact manifold with
boundary there are a number of different natural algebras of
pseudodifferential operators which generalize the `standard' algebra of
pseudodifferential operators on a compact manifold without
boundary. Amongst these are the calculus of b-pseudodifferential operators
\cite{Melrose25} (b=boundary), \cite{Melrose42}, the scattering calculus
\cite{Melrose43} and the uniformly degenerate calculus (or zero)
\cite{Mazzeo2} and \cite{Mazzeo-Melrose1}.  The distinction between the
terms `calculus' and `algebra' is not great here. The former is preferred
because all of the algebras we discuss have natural, and useful, extensions
to somewhat larger spaces of operators in which not every pair of elements
can be composed.  If the manifold has more structure, for example if its
boundary admits a fibration, then there are other possibilities, such as
the edge calculus \cite{Mazzeo4} which interpolates between the b and
uniformly degenerate calculi. In this paper we shall discuss another
algebra of this general type; it is associated to a fibration of the
boundary and a choice of boundary defining function up to second order at
the boundary, or more precisely to a trivialization of the conormal bundle
to the boundary over each fibre. The extreme cases, in terms of the fibre
dimension of the fibration, of this algebra correspond to the `cusp'
algebra, of operators naturally associated to (finite volume) hyperbolic
cusps, and the scattering algebra, of operators associated to Euclidean
scattering theory.

The purpose of this paper is to give a concise yet complete treatment of
this `fibred-cusp' algebra, along with a few of the most basic
consequences. More sophisticated applications will be taken up elsewhere.
In this introduction we shall give an outline of some of the salient
features of the algebra which will be proved in full later in the paper.

Let $X$ be a compact \ci\ manifold with boundary and suppose that the
boundary has a smooth fibration 
\begin{equation}
\phi :\pa X\longrightarrow Y,
\label{I.1}
\end{equation}
where $Y$ is the space of fibres. Suppose also that $x\in\CI(X)$ is
a choice of boundary defining function, \ie $x\ge0,$ $\pa X=\{x=0\}$ and
$dx\not=0$ at $\pa X.$ In particular, $x$ fixes a trivialization of the
conormal bundle to the boundary. Associated with this structure is the
space of {\em fibred cusp} vector fields 
\begin{equation}
\begin{aligned}
\mathcal{V}_{\Phi}(X)=\big\{V\in\CI(X;TX);&Vx\in x^2\CI(X)\Mand V_p\\
&\text{is tangent to }\phi ^{-1}(\phi (p))\ \forall\ p\in\pa X\big\}.
\end{aligned}
\label{I.2}
\end{equation}

As shown below, $\mathcal{V}_{\Phi}(X)$ is a Lie algebra and $\CI(X)$ module
which is projective in the sense that there is a \ci\ vector bundle $\FT X$
over $X$ with natural vector bundle map $\iota _{\Phi}:\FT X\longrightarrow
TX,$ which is an isomorphism over $X^{\circ}=X\setminus\pa X,$ and is such
that 
\begin{equation*}
\CI(X;\FT X)=\iota _{\Phi}\circ\mathcal{V}_{\Phi}(X).
\label{I.3}\end{equation*}
That is, $\mathcal{V}_{\Phi}(X)$ can be naturally identified with $\CI(X;\FT
X).$ The identifier `$\Phi$' will be used to denote objects which are
naturally associated to $\mathcal{V}_{\Phi}(X)$. Note that 
$\mathcal{V}_{\Phi}(X)$ determines the map $\phi$, but does not
completely determine the defining function $x$.

There are two extreme cases to keep in mind as a guide to this 
discussion, occurring when $\phi$ is one of the `trivial' (or
universal) fibrations. The first is when $Y=\{\text{pt}\}$ and the
second when $Y=\pa X$. In the former case, $\mathcal{V}_{\Phi}(X)$
determines, and is determined by, the defining function $x$ up to the
equivalence $x'\sim x$ if $x'=cx+x^2g,$ where $c>0$ is constant and
$g\in\CI(X).$ This will be called the {\em cusp algebra.} In the
latter case, the Lie algebra is independent of the choice of $x$ and 
is called the {\em scattering algebra}. The algebra of pseudodifferential 
operators associated to it is discussed in \cite{Melrose43} and 
\cite{Melrose-Zworski1} and in local form on $\RR^n$ goes back at least to
Shubin \cite{Shubin1}. When $X$ is the upper half-sphere, the
interior of which may be identified with $\RR^n$ via
stereographic compactification $\RR^n\hookrightarrow \SS^n_+$,
the scattering algebra is generated by the translation-invariant vector 
fields.

Since $\mathcal{V}_{\Phi}(X)$ is a Lie algebra and $\CI(X)$ module it is
natural to consider the enveloping algebra, $\DiffF*(X),$ consisting of
those operators on $\CI(X)$ which can be written as finite sums of products
of elements of $\mathcal{V}_{\Phi}(X)$ and $\CI(X).$ It is filtered by the
subspaces $\DiffF k(X)$ which have elements expressible as sums of products
involving at most $k$ factors from $\mathcal{V}_{\Phi}(X).$ Let $\FT^*X$ be the
dual bundle to $\FT X$ and let $P^k(\FT^*X)\subset\CI(\FT^*X)$ be the space
of functions which are homogeneous polynomials of degree $k$ on the
fibres. The principal symbol map extends from the interior to $\sigma
_{\Phi,k}:\DiffF k(X)\longrightarrow P^k(\FT^*X)$. This map is multiplicative
and gives a short exact sequence delineating the filtration
\begin{equation}
0\longrightarrow \DiffF{k-1}(X)\hookrightarrow \DiffF k(X)\overset{\sigma
_{\Phi,k}}\longrightarrow P^k(\FT^*X)\longrightarrow 0.
\label{I.4}
\end{equation}

We {\em microlocalize} this algebra of differential operators to obtain the
filtered algebra of fibred-cusp, or $\Phi$-, pseudodifferential operators 
\begin{equation*}
\DiffF k(X)\subset\PsiF k(X)
\label{I.5}\end{equation*}
where $\PsiF m(X)$ is defined for each $m\in\RR.$ Again there is a
multiplicative symbol map delineating the filtration 
\begin{equation}
0\longrightarrow\PsiF{m-1}(X)\hookrightarrow\PsiF m(X)
\overset{\sigma_{\Phi,m}}\longrightarrow S^m(\FT^*X)/S^{m-1}(\FT^*X)
\longrightarrow 0
\label{I.6}
\end{equation}
where $S^m(E),$ for any vector bundle $E,$ is the space of symbols of order
$m.$ The construction of $\PsiF m(X)$ is effected geometrically. More
specifically, these spaces of operators are characterized by 
the regularity properties of their Schwartz' kernels. These, in turn,
are defined as conormal distributions on a space, $\XF2,$ which is
a resolution of $X^2.$ This resolution is obtained from the
ordinary `double space' through a sequence of blow-ups. 
One of the main facts about $\PsiF{*}(X)$, that it is closed
under composition, is proved using a resolution $\XF3$ of the
ordinary triple space $X^3$, as we shall explain later.

Whether a particular element in $\PsiF{*}(X)$ acts as a Fredholm operator, 
say on $L^2$, is no longer determined solely by the invertibility of its 
image under the symbol map \eqref{I.4} or \eqref{I.6}.
In fact, there is a second symbol map, the range of which is in general 
no longer a commutative algebra. To introduce this {\em normal operator}, we 
first describe the space of operators in which it lies.

If $F$ is any compact manifold without boundary and $W$ is a real vector
space, then the space $\Psi^m(F\times W)$ of all pseudodifferential operators 
on the \ci\ manifold  $F\times W$ is well defined. This is not
an algebra because we have imposed no growth restrictions on the kernels. A
special subclass consists of those 
elements which are invariant under translation in $W,$
and therefore loosely speaking act by convolution in the $W$ factor and as
ordinary pseudodifferential operators in $F$. Now consider
\begin{equation}
\Psis Wm(F)\subset\Psi^m(F\times W)
\label{I.7}
\end{equation}
consisting of those translation invariant operators with convolution 
kernels on $F^2\times W$ which are rapidly decreasing with all derivatives 
at infinity. These spaces form a filtered algebra in the usual way and 
we call them the `$W$-suspended pseudodifferential operators on $F$', even 
though they act on functions on $F\times W.$ They are invariant under 
diffeomorphisms of $F$ and linear transformation of $W.$ This means
that we can define $\Psisf Wm(X';W)$, where $\phi:X'\longrightarrow Y$ is any 
fibration, $W\longrightarrow Y$ a vector bundle, and $G=X'\times_{Y}W$ 
the fibre product, where elements are defined as in \eqref{I.7} on the 
fibres of $G$ and depend smoothly on the base variable in $Y.$

If $\phi:\pa X\longrightarrow Y$ is the fibration \eqref{I.1}, and 
$\iota _{\Phi}:\FT X\longrightarrow TX$ is the natural inclusion
map, set 
\begin{equation}
\FN_p\pa X=\left\{v\in\FT_pX,\ p\in\pa X;\iota _{\Phi}(v)=0\right\}.
\label{I.8}
\end{equation}
Although this is defined as a bundle over $\pa X$, in fact it is the 
lift to $\pa X$ of a bundle, $\FN Y,$ over $Y,$
\begin{equation}
\FN\pa X=\phi ^*(\FN\pbX),
\label{I.14}
\end{equation}
and hence is of the form just described. The {\em normal homomorphism},
which we will define later, takes values in the corresponding space of
suspended operators, and there is a multiplicative short exact sequence
\begin{equation*}
0\longrightarrow x\PsiF m(X)\hookrightarrow \PsiF
m(X)\overset{\NF}{\longrightarrow }\Psisf {\FN Y}m(\pa X)\longrightarrow 0.
\label{I.9}\end{equation*}

The symbol and normal operator together are sufficient to capture the 
Fredholm property for these differential or pseudodifferential operators.

\begin{theorem}\label{I.10} An element $P\in\PsiF0(X)$ is Fredholm as an
operator on $L^2(X)$ if and only if it is fully elliptic in the sense
that its symbol $\sigma_{\Phi,0}$ is invertible and in addition its 
normal operator $\NF(P)$ is invertible as an element of $\Psisf {\FN
Y}m(\pa X).$
\end{theorem}

We will state and prove a more general result for pseudodifferential
operators of any order acting on sections of a vector bundle. This raises
the following fundamental

\begin{prob}\label{I.11} Find an explicit index formula for fully elliptic
$\Phi$-pseudodifferential operators in terms of the symbol and normal operator.
\end{prob}

This has been done in full generality in only one case, where $Y=\pa X,$
\ie for the scattering calculus. This is discussed briefly in
\cite{Melrose44}, where it is reduced to the Atiyah-Singer theorem. In the
other extreme case, where $Y=\{\text{pt}\},$ the calculus is essentially 
that of manifolds with cylindrical ends. The index theorem in this
setting for Dirac operators is that of Atiyah, Patodi and Singer
\cite{Atiyah-Patodi-Singer1}.  There is a somewhat non-explicit index
formula for general fully elliptic pseudodifferential operators here due to 
Piazza \cite{Piazza1}. In \cite{Melrose46} a definition of the eta
invariant in this context is given, and \cite{Melrose-Nistor2} contains
an index formula in terms of it.

Beyond these index questions, another reason for developing these 
calculi of operators is to analyze the regularity of solutions 
to related differential equations. We formalize this process using
the notion of a {\em wavefront set}, which is defined by microlocal 
invertibility properties of $\Phi$-pseudodifferential operators. In the
analytic category the wavefront set (singular spectrum) was introduced by
Sato, see \cite{Kashiwara-Kawai-Sato1}; in the \ci\ category it is due to
H\"ormander \cite{Hormander1}.

To describe this consider again the structure bundle $\FT X$ and 
its dual $\FT^*X.$ The stereographic compactification of a vector space 
to a ball, or half-sphere, is linearly covariant, and so we can
define the fibrewise compactification any vector bundle. Since $X$ is a
manifold with boundary the compactification $\cFTs X$ is a manifold with
corners up to codimension two. The restriction to the boundary of the
bundle $\FT^*X$ has as quotient $\FN^*\pa X$ which is, as noted above,
naturally the lift of a bundle $\FN^*Y$ over the base $Y.$ This is the
parameter space for the normal operator. The disjoint union of the part `at
infinity' of the bundle $\cFTs X$ and the compactification, $\cFN^*Y,$
\begin{equation*}
\CF=\FS^*X\sqcup\cFN^*Y,
\label{I.12}\end{equation*}
is the carrier of the $\Phi$-wavefront set 
\begin{equation*}
\WFF(u)=\WFF^{\sigma}\sqcup\WFF^{\pa}\subset\CF.
\label{I.13}\end{equation*}
It has properties and utility similar to the usual wavefront set.

The authors thank Andr\'as Vasy for helpful comments.

\paperbody
\section[Fibred cusp algebras]{Fibred cusp algebras \label{S.FCA}}

We begin our more detailed discussion by analyzing the space of vector fields 
defined by \eqref{I.2}. Thus, $X$ is a compact \ci\ manifold with boundary 
and as in \eqref{I.1}, $\phi$ is a fibration of the boundary. If the boundary 
is not connected we denote by $M_1(X)$ the set of boundary components. Then 
each boundary hypersurface $H\in M_1(X)$ has a specified fibration 
$\phi_H:H\longrightarrow Y_H.$ There need be no relationship between these
fibrations. For the most part we shall simplify the discussion by supposing
that $\pa X$ is connected, but when confusion might arise in the general case
we make a precise statement.

In addition to the fibration, we also suppose that a boundary defining 
function $x\in\CI(X)$ is given. As will be discussed shortly, the 
structure we describe does not depend on all the information in $x.$ 
Consider $\mathcal{V}_\Phi (X)$ defined by \eqref{I.2} which should now 
be written more carefully as
\begin{equation*}
\begin{aligned}
\mathcal{V}_{\Phi}(X)=\big\{V\in&\, \CI(X;TX);Vx\in x^2\CI(X)\Mand V_p\\
&\text{is tangent to }\phi_H ^{-1}(\phi (p))\ \forall\ p\in H,\ \forall\ H\in
M_1(X)\big\}.
\end{aligned}
\end{equation*}

\begin{lemma}\label{LCD} Suppose $p\in \pa X$ and $y_1,\dots, y_{\ell}$ are
local coordinates in $Y$ near $\phi(p).$ Let $\widetilde y_1,\dots, \widetilde
y_\ell\in\CI(X)$ be functions satisfying $\widetilde
y_j=\phi^*(y_j)$ on $\pa X$ near $p$ and choose $k=n-\ell-1$
functions $z_1,\dots, z_k$ such that $x,\widetilde{y_j}, z_i$ give
local coordinates in $X$. Then near $p,$ $\mathcal{V}_\Phi (X)$ is spanned by
\begin{equation}
x^2\frac{\pa}{\pa x}, \ x\frac{\pa}{\pa \tilde y_j}, \ \frac{\pa}{\pa z_i}.
\label{FCA.1}
\end{equation}
\end{lemma}

\begin{proof} Since the differentials of $x$ and the $\widetilde{y}_j$ must be
independent at $p$ there do indeed exists functions $z_i$ completing them
to a coordinate system. A general vector field on $X$ is locally
\begin{equation*}
V=a\frac{\pa}{\pa x}+\sum\limits_{j=1}^\ell b_j\frac{\pa}{\pa\tilde y_j}+
\sum\limits_{i=1}^k c_i \frac{\pa}{\pa z_i}
\end{equation*}
for \ci\ coefficients $a,b_j, c_i.$ Then $Vx=a,$ so the first
condition on $V$ in \eqref{I.2} is that $a=O(x^2),$ \ie $a=x^2 a'$
with $a'$ \ci\ near $p.$ Locally the fibres of $\phi$ are the
surfaces $\tilde y=$const, in $x=0.$ Thus if $V\in\mathcal{V}_\Phi (X)$
then $b_j=xb_j'.$ This shows that the elements in \eqref{FCA.1} span
$\mathcal{V}_\Phi (X)$ locally over $\CI(X).$
\end{proof}

Lemma~\ref{LCD} actually shows that $\mathcal{V}_\Phi (X)$ is projective,
and this means that we can interpret this space of vector fields as the
full set of sections of some vector bundle. For any $p\in X$ let 
$\mathcal{I}_p (X)\subset \CI(X)$ be the ideal of functions vanishing 
at $p.$  Then denote by $\mathcal{I}_p\cdot\mathcal{V}_\Phi (X)\subset
\mathcal{V}_\Phi (X)$ the finite linear span of products, $aV$, for
$a\in\mathcal{I}_p(X)$ and $V\in\mathcal{V}_\Phi (X),$ and set
\begin{equation*}
\FT_pX=\mathcal{V}_\Phi (X)\big/\mathcal{I}_p\cdot\mathcal{V}_\Phi(X).
\end{equation*}

\begin{lemma}\label{TSB} For each $p\in X,$ $\FT_pX$ is a vector space
of dimension $\dim X,$ and the disjoint union
\begin{equation*}
\FT X=\bigsqcup\limits_{p\in X}\FT_pX
\end{equation*}
has a natural structure as a smooth vector bundle over $X.$ 
There is a natural linear map $\iota_p:\FT_pX\longrightarrow T_pX$ 
which is an isomorphism when $p\in X^{\circ}=X\backslash \pa X;$ these maps
define a smooth bundle map $\iota:\FT X\longrightarrow TX$ with 
the property that for every $\tilde V \in \CI(X;\FT X)$ there is a unique 
$V\in\mathcal{V}_\Phi (X)\subset \CI(X;TX)$ such that
\begin{equation*}
\iota_p\tilde V_p=V_p \ \forall \ p\in X^{\circ}.
\end{equation*}
Conversely, each $V\in\mathcal{V}_\Phi (X)$ defines a section
$\tilde V\in\CI(X;\FT X).$ 
\end{lemma}

\begin{proof} Over the interior of $X$ the elements of $\mathcal{V}_\Phi (X)$ 
are unconstrained, and so $\FT_pX\equiv T_pX$ for $p\in X^{\circ}.$ We write 
this identification as $\iota _p:\FT_pX\longrightarrow T_pX.$  Near a 
boundary point $p$ we have shown that $V\in\mathcal{V}_\Phi (X)$ has a unique 
smooth decomposition in terms of the vector fields \eqref{FCA.1}. Thus
$V'\in\mathcal{I}_p\cdot \mathcal{V}_\Phi (X)$ if and only if its
decomposition has coefficients vanishing at $p.$ This means that 
(the residue classes of) $x^2\frac{\pa}{\pa x},$ $x\frac{\pa}{\pa\tilde y_j}$ 
and $\frac{\pa}{\pa z_i}$ give a basis of $\FT_pX$, and therefore
this vector space has dimension $\dim X.$  In fact, these sections clearly 
give $\FT X$ the structure of a vector bundle near $p,$ where any smooth 
section is (locally) given by an element of $\mathcal{V}_\Phi (X)$ and
conversely. It remains only to show that this vector bundle structure is
independent of the choice of the local coordinates. This follows simply by
inserting the change of coordinate formula for vector fields into the basis
\eqref{FCA.1}.
\end{proof}

Generally we shall ignore the map $\iota$ and identify $\mathcal{V}_\Phi 
(X)$ with $\CI(X;\FT X)$ as this lemma permits us to do. As noted in 
the introduction

\begin{lemma}\label{LAP} The space $\mathcal{V}_\Phi (X)$ is a Lie subalgebra of
$\CI(X;TX).$
\end{lemma}

\begin{proof} If $V,$ $W\in\mathcal{V}_\Phi(X)$, then by definition 
they are tangent to the fibres of $\phi$ in $\pa X$. Because 
tangency to a submanifold persists for commutators, $[V,W]$ also
has this property. Similarly, since $Vx = x^2 a$ and $Wx = x^2 b$
for some functions $a,b \in \CI(X)$, 
\begin{multline*}
[V,W]x=V(Wx)-W(Vx)=V(x^2 a)-W(x^2b)\\
=x^2(Va-Wb) + 2x(a-b)(V-W)(x) = x^2(Va-Wb) + 2x^3(a-b)^2.
\end{multline*}
\end{proof}

As noted above, the algebra $\mathcal{V}_{\Phi}(X)$ is determined by $\phi$
and the choice of a boundary defining function. Conversely,
$\mathcal{V}_{\Phi}$ determines $\phi$ but it does not completely determine
$x.$ In fact two boundary defining functions $x$ and $x'$ determine the
same Lie algebra $\mathcal{V}_{\Phi}(X),$ and hence the same `boundary
structure' relative to $\phi,$ if and only if $x'=\alpha x,$ with $\alpha
|_{\pa X}\in\phi ^*\CI(Y).$ Thus if we let $\CI_\phi(X)\subset\CI(X)$ be
the space of smooth functions on $X$ which are constant on each leaf of
$\phi$ at the boundary then this means $x'\in x\CI_{\phi}(X).$

The Lie algebra $\mathcal{V}_{\Phi}(X)$ has a natural ideal, consisting of
those elements which vanish at the boundary as vector fields in the usual
sense. In terms of the basis \eqref{FCA.1}, it is spanned by $x^2\pa_x$
and $x\pa_{\tilde y_j}$ and $x\pa_{z_i}$ near each boundary point and is
unconstrained in the interior. Over the boundary it spans the subbundle
\eqref{I.8}. This ideal is the span over $\CI(X)$ of a smaller subalgebra  
\begin{equation}
\mathcal{W}_{\Phi}(X)=\{V\in\mathcal{V}_{\Phi}(X);V\in x\CI(X,TX)\Mand
Vx\in x\CI_{\phi}(X)\}.
\label{FCA.100}\end{equation}
This latter condition is clearly independent of the choice of $x$ defining
$\mathcal{V}_{\Phi}(X),$ \ie $\mathcal{W}_{\Phi}(X)$ \emph{is} an
invariantly defined subspace of the latter. Since $\mathcal{W}_{\Phi}(X)$
is a $\CI_{\phi}(X)$-module (and not a $\CI(X)$ module) the subbundle of
$\FT_{\pa X}X$ it defines is naturally the lift of a bundle from $\pbX.$
This is the bundle $\FN\pbX$ in \eqref{I.14}. 

There is a direct representation of the fibre $\FT^*_pX$, $p\in\pa X$, of
the dual bundle which is useful later. For $p\in\pa X$ let
$I_p(X)\subset\CI_\phi(X)$ be the ideal of functions vanishing on the fibre
through $p$ and $J_p(X)\subset I_p(X)$ the smaller ideal of functions with
restriction to $\phi^{-1}(p)$ vanishing to second order at $p.$ Note that
$I_p(X)$ only depends on the fibre through $p$ but $J_p(X)$ depends also on
the location of $p$ within this fibre. If $x$ is an admissible boundary
defining function then there is a canonical isomorphism
\begin{equation*}
\FT^*_pX\equiv x^{-1}\CI_\phi(X)/x^{-1}J_p(X)
\end{equation*}
given by applying $V\in\mathcal{V}_{\Phi}(X)$ and evaluating at $p;$ this
follows from \eqref{FCA.1}.

Let $\DiffF m(X)$ be the space of operators on $\CI(X)$ generated by 
$\CI(X)$ and products of up to $m$ elements of $\mathcal{V}_\Phi (X)$.
The local structure of these operators is easy to determine.

\begin{lemma}\label{SEA} In the local coordinates near a boundary point
described in Lemma~\ref{LCD}, any $P\in \DiffF m(X)$ may be written 
\begin{equation}
P=\sum\limits_{|\alpha |+|\beta |+q\le m} p_{\alpha ,\beta ,q} (x,y,z)(x^2
D_x)^q(xD_y)^\beta D^\alpha _z, \qquad D_t=\frac{1}{i} \frac{\pa}{\pa t}.
\label{FCA.6}
\end{equation}
Conversely, if $P\in\Diff m(X)$ and this holds in a neighbourhood of 
each boundary point, then $P\in\DiffF m(X).$
\end{lemma}

\begin{proof} This follows by induction on $m.$ Certainly \eqref{FCA.6}
holds when $m=0.$ In general, $\DiffF{m+1}(X)$ is the span of
$\mathcal{V}_{\Phi}(X)\cdot\DiffF m(X)$ and $\DiffF m(X).$
Using the local representation of $\mathcal{V}_{\Phi}(X)$ given by
\eqref{FCA.1} and the representation \eqref{FCA.6} for $\DiffF m(X)$, the
same result follows directly for $\DiffF{m+1}(X).$
\end{proof}

Note that the order of the various factors in \eqref{FCA.6} is immaterial, 
because changing it would just change the coefficients slightly. 
The main properties of this space of $\Phi$-differential operators will
be discussed in more detail once we have defined the 
space of $\Phi$-pseudodifferential operators. 

\section[$\Phi$-pseudodifferential operators]{$\Phi$-pseudodifferential 
operators\label{S.FDO}}

We now turn to the definition of the `small' calculus of
$\Phi$-pseudodifferential operator. These can be thought of as `symbolic'
functions of the vector fields in $\mathcal{V}_\Phi (X)$ in the same sense 
that, by \eqref{FCA.6}, the $\Phi$-differential operators are polynomial
functions in these vector fields. Our definition of this calculus
is quite geometric; this has the virtue that many of
the main properties we need to develop, in particular
the fact that this space of operators is closed under
composition, may be proved directly and also quite
geometrically. 

Following a general `microlocalization' principle for algebras 
of this type, the operators in $\mathcal{V}_\Phi (X)$ will be 
characterized by the lifts of their Schwartz kernels from
$X^2$ to a space $X^2_\Phi$ which is obtained by a resolution process, 
more specifically by blowing up a sequence of p-submanifolds in $X^2$. Here
the p-submanifolds (for `product') are those around which the manifold with
corners has a product decomposition, they may be thought of as properly
embedded. The point of this geometric resolution is that it encodes the
approximate local homogeneities of $\Phi$-differential operators,
and so it is natural to define the $\Phi$-pseudodifferential operators
by requiring that their Schwartz kernels also have the same
approximate local homogeneities, \ie lift to well-behaved
distributions on $X^2_\Phi$. We refer to 
\cite{Epstein-Melrose-Mendoza1}, \cite{Melrose40}, \cite{Mazzeo-Melrose1} and
\cite{Hassell-Mazzeo-Melrose1} for a discussion of the process
of blowing up a p-submanifold in a manifold with corners.

As already noted, all our constructions proceed independently
at each boundary hypersurface of $X$, and so it is sufficient to 
suppose that $\pa X \equiv H$ is connected.

The Schwartz kernel of any operator on $\CI(X)$ 
is a distribution on $X^2$. Of course, we are particularly interested 
in the behaviour of these operators, and hence kernels, near the boundary.
We use the notation 
\begin{equation*}
L(H)=H\times X, \ R(H)=X\times H
\end{equation*}
or simply $L,R$ when $H$ is understood. For any manifold with corners
$Z$, let $M_k(Z)$ denote set of boundary components of codimension $k$.
In particular $\{L,R \} = M_1(X^2)$. Because $X$ is a manifold with
boundary, $X^2$ has boundary components only up to codimension two.
Amongst these, only the faces $B(H) \in M_2(X^2)$, 
\begin{equation*}
B=B(H)=H\times H\subset X^2, 
\end{equation*}
which are the ones intersecting the diagonal, are of interest to us.
The other manifold of primary importance in this discussion is the diagonal
\begin{equation*}
D=\{(z,z)\in X^2\}.
\end{equation*}

An important feature of this geometry is that these submanifolds 
do not intersect normally. We resolve this by blowing up $B$ to get
the b-double space
\begin{equation}
X^2_b=[X^2;B]; \ \beta ^2_b:\ X^2_b\longrightarrow X^2.
\label{FDO.4}
\end{equation}
This compact manifold with corners is obtained by taking the disjoint 
union of $X^2 \setminus B$ and the inward-pointing spherical normal
bundle of $B$ and endowing this set with the (unique) 
minimal \ci\ structure for which smooth functions on $X^2$ 
and polar coordinates in $X^2$ around $B$ all lift to be
smooth. 

We next describe the lifts of the submanifolds of $X^2;$ we 
shall use the same letters to denote the lifts but add a 
subscript `b' when necessary to distinguish between a manifold 
and its lift.

First, the `front face' of $X^2_b$, which is produced by the blow-up 
of $B$ is denoted 
\begin{equation*}
B_b=(\beta ^2_b)^{-1}(B)=SN^+B.
\label{FDO.5}
\end{equation*}
By definition it is a quarter circle bundle over $B$. In fact, since 
the fibred-cusp structure specifies the $1$-jet of a defining function $x$ 
for the boundary $H$ of $X,$ this bundle is naturally trivial over $B$
\begin{equation}
B_b=B\times [-1,1]_s.
\label{FDO.6}
\end{equation}
To see this, note that if $x$ and $x'$ are the lifts to $X^2$ of the given 
boundary defining function from the left and right factors of $X$,
respectively, then $NB$ is spanned by $\pa_x$ and $\pa_{x'}$. The interior 
normal bundle $N^+B$ is therefore $\{(p,a\pa_x+a'\pa_{x'}); a,a'\ge 0\}$ 
and it is then easy to check that \eqref{FDO.6} follows, if we use
$s=\frac{a-a'}{a+a'}.$

Next consider $L$ and $R;$ the inverse images of these boundary faces under
$\beta^2_b$ contain $B_b.$ We define instead their {\em lifts} to exclude
the interior of the front face:
\begin{equation*}
\begin{gathered}
L_b=\clos\bigg((\beta ^2_b)^{-1}(L)\big\backslash
B_b\bigg)=\clos\bigg((\beta ^2_b)^{-1}(L\big\backslash B)\bigg)\\
R_b=\clos\bigg((\beta ^2_b)^{-1}(R)\big\backslash
B_b\bigg)=\clos\bigg((\beta ^2_b)^{-1}(R\big\backslash B)\bigg)
\end{gathered}
\label{FDO.7}
\end{equation*}
where $\clos$ denotes the closure. These are boundary hypersurfaces of
$X^2_b$ and all of the boundary hypersurfaces of $X^2_b$ have been
enumerated, so that 
\begin{equation*}
M_1(X^2_b)=\bigcup\limits_{H\in M_1(X)}\{L_b(H), R_b(H), B_b(H)\}.
\label{FDO.8}
\end{equation*}

We also define the lifted diagonal
\begin{equation*}
D_b=\clos\left((\beta ^2_b)^{-1}(D\big\backslash D \cap B)\right)
\subset X^2_b.
\label{FDO.9}
\end{equation*}
As noted earlier, the diagonal itself does not intersect the
boundary normally. However, $D_b$ is a closed embedded p-submanifold,
and the only boundary hypersurfaces it meets are the diagonal
front faces $B_b(H).$ In this sense the blow-up of $X^2$ to $X^2_b$ 
resolves the `geometry' on $X^2$ consisting of the boundary 
faces and the diagonal. This blow-up is the basis for the direct definition
of the b-calculus, see \cite{Melrose42}.

There are however further degeneracies, associated to the fibred-cusp
algebra, which need to be resolved. These occur along the fibre diagonal
of the front face $B=H\times H$, given by 
\begin{equation*}
\{(h,h')\in B; \ \phi(h)=\phi(h') \text{ in } Y\}.
\end{equation*}
Using the product decomposition \eqref{FDO.6} this lifts to the
submanifold 
\begin{equation}
\Phi =\Phi (H) \equiv \{(h,h',0)\in B_b=B\times[-1,1]_s;\phi(h)=
\phi(h')\},
\label{FDO.10}
\end{equation}
which is an embedded, closed submanifold in the interior of
$B_b,$ hence is a p-submanifold of $X^2_b.$ For any manifold with corners
$Z$ we denote by $\mathcal{V}_b(Z)$ the Lie algebra of smooth vector fields
which are tangent to each of the boundary faces.

\begin{prop}\label{CLA} The Lie algebra $\mathcal{V}_b(X)$ lifts to
$X^2_b$ from either factor to a Lie subalgebra of $\mathcal{V}_b(X^2_b)$
transversal to $D_b.$ In each case the elements of the lift of
$\mathcal{V}_\Phi(X)$ constitute the subset of the lift of
$\mathcal{V}_b(X)$ consisting of those vector fields which are tangent to
$\Phi.$
\end{prop}

\begin{proof} Using local coordinates $x,$ $y,$ $z,$ $x',$ $y',$ $z',$
as in Lemma~\ref{LCD}, on both the left and right factors of $X$, 
gives coordinates 
\begin{equation}
x',s=\frac{x-x'}{x+x'}, \ \tilde y, \tilde y', z,z' \Min X^2_b,
\label{FDO.11}
\end{equation}
on $X^2_b$ valid near $s=0$. The vector fields in \eqref{FCA.1} lift to
\begin{equation*}
\frac{x'}{2}(1+s)\pa_s,\frac{2x'}{1-s}\frac{\pa}{\pa\tilde y_j},
\frac{\pa}{\pa z_j}, 
\label{FDO.12}\end{equation*}
and these are clearly tangent to $\Phi =\{x'=0,\ s=0, \tilde y=\tilde y'\}.$
On the other hand, the basic generating set of vector fields on 
$\mathcal{V}_b(X)$ is
\begin{equation}
x\pa_x, \ \pa_{\tilde y_{j}}, \ \pa_{z_i};
\label{FDO.13}
\end{equation}
these lift to
\begin{equation*}
\frac 12(1-s^2)\pa_s, \ \pa_{\tilde y_j},\ \pa_{z_{j}}
\label{FDO.14}\end{equation*}
which are clearly transversal to $D_b=\{s=0,\tilde y=\tilde y', z=z'\}.$ If
such a lift is tangent to $\Phi$ then it is easily seen to be the lift of
an element of $\mathcal{V}_\Phi(X).$
\end{proof}

It follows from \eqref{FDO.11} that $\Phi$ is the flow out 
of $\pa D_b$ under the lifts of $\mathcal{V}_\Phi (X)$ from the left 
and right factors. It is therefore the minimal submanifold to which these 
lifted vector fields are tangent.


In the second (and final) stage of the fibred boundary blow-up we 
define 
\begin{equation}
X^2_\Phi =[X^2_b; \Phi ], \qquad \beta_{\Phi-b}, X^2_\Phi \longrightarrow X^2_b.
\label{FDO.15}
\end{equation}
There is also a full blow-down map
\begin{equation*}
\beta_{\Phi} = \beta_{\Phi-b}\circ \beta_b: X^2_\Phi \longrightarrow X^2.
\label{FDO.155}
\end{equation*}

\begin{lemma}\label{TLL} The Lie algebra $\mathcal{V}_\Phi (X)$ lifts,
from either the left or right factor of $X$ to a Lie subalgebra 
of $\mathcal{V}_b(X^2_\Phi ).$  The diagonal $D$ lifts to a closed embedded 
p-submanifold 
\begin{equation}
D_\Phi =\clos\beta^{-1}_\Phi (D\cap (X^2)^{\circ})
\label{FDO.16}
\end{equation}
and the lifted algebra is transversal to the lifted  diagonal.
\end{lemma}

\begin{proof} These statements are all trivial in the interior of $X^2_\Phi,$
which is diffeomorphic to the interior of $X^2$. Since all constructions
are local near fibres of $\phi$ it suffices to consider the model
product fibred-cusp structure 
\begin{equation*}
X=[0,1)_x\times Y_y\times F_z.\label{FDO.17}\end{equation*}
Near its front face, $X^2_b$ is also a product
\begin{equation}
X^2_b\simeq [-1,1]_s \times[0,1)_{x+x'}\times Y^2\times F^2.
\label{FDO.18}
\end{equation}
The boundary fibre diagonal is the submanifold
\begin{equation*}
\Phi =\{0\}\times\{0\}\times\Diag_Y\times F^2.
\label{FDO.19}\end{equation*}

The second blow-up occurs only in the first three factors of 
\eqref{FDO.18}. This effectively reduces the problem to the 
case $F=\{\operatorname{pt}\}.$ In the first blow-up, 
\eqref{FDO.4}, if $r=x+x'$ then
\begin{equation*}
x=\frac 12(s-1)r, \ x'=\frac 12 (1-s)r, \ s=\frac{x-x'}{x+x'}.
\label{FDO20}\end{equation*}
As noted in Proposition~\ref{CLA}, the basis fields \eqref{FCA.1} 
lift smoothly to $X^2_b$ as
\begin{equation*}
\begin{gathered}
x^2\pa_x\longmapsto \frac{1}{4}(1+s)^2 r^2\pa_r+\frac{1}{4}(1+s)^2(1-s)
r\pa_s\\
x \pa_{\tilde y_j}\longmapsto \frac 12 (1+s)r\pa_{\tilde y_j}.
\end{gathered}
\label{FDO.21}
\end{equation*}
Since the lifts of smooth coefficients are smooth, we can instead
consider near $\Phi$ the simpler basis 
\begin{equation}
r(r\pa_r+(1-s)\pa_s), \ r\pa_{\tilde y_j},
\label{FDO.22}
\end{equation}
which also spans the lift of $\mathcal{V}_\Phi (X)$ over $\CI(X^2_b)$.

Since these vector fields are tangent to $\Phi =\{r=0, s=0, y=\tilde y\},$  
they lift smoothly under the blow-up of $\Phi.$ Near the lifted diagonal 
the variables $r,$ $S=s/r,$ $Y_j=(\tilde y_j-\tilde y'_{j})/r$ give local 
coordinates and in terms of these the vector fields in \eqref{FDO.22} become 
\begin{equation*}
r^2\pa_r-rS\pa_S-rY\cdot \pa_Y+(1-r S)\pa_S, \ \pa_{Y_{j}}.
\label{FDO.23}\end{equation*}
Since the lifted diagonal itself is $\{S=0, \ Y=0\}$, this shows that 
the lift of $\mathcal{V}_\Phi (X)$ from the left factor (and hence by
symmetry also from the right factor) is transversal to the lifted 
diagonal $D_\Phi.$
\end{proof}

Notice that $D_\Phi \simeq X,$ with the diffeomorphism given by 
$\pi\circ\beta _\Phi^2$. The transversality in Lemma~\ref{TLL} shows that
there are natural isomorphisms from the normal and conormal bundles
of $D_\Phi$ in $X^2_\Phi$ to the $\Phi$-tangent and cotangent
bundles of $X$  covering this identification,
\begin{equation*}
\FT X\simeq ND_\Phi, \ \FT^*X\simeq N^*D_\Phi. \label{FDO.24}
\end{equation*}

Now that $X^2_\Phi$ has been defined, we consider the structure of the
lifts of the Schwartz kernels of $\Phi$-differential operators
to this space. The fundamental case to understand is the identity
operator. In local coordinates its Schwartz kernel (on $X^2$) is 
\begin{equation*}
\begin{aligned}
K_{\Id}=\delta (x-x')\delta (y-y')\delta (z-z')\, dx'dy'dz'\\
\Longrightarrow u(x,y,z)=\int K_{\Id} u(x', y', z').
\end{aligned}
\label{FDO.25}
\end{equation*}
From now on we identify this distribution with the operator $K_{\Id}=\Id.$
It is more convenient to write it in terms of a `$\Phi$-density' 
$(x')^{-(\ell+2)}dx'dy' dz'$ where $\ell=\dim Y:$
\begin{equation}
\begin{gathered}
\Id=(x')^{\ell+2}\delta (x-x')\delta (y-y')\delta (z-z')\nu
'_\Phi \\ 
\nu '_\Phi =(x')^{-\ell-2}dx' dy' dz'.
\end{gathered}
\label{FDO.26}
\end{equation}

Consider its behaviour when lifted from the interior of $X^2$ to $X^2_b$.
It is supported on the diagonal $D,$ so it suffices to consider a
neighbourhood of $\{s=0\}$ in \eqref{FDO.16}.  Since $x-x'=rs$ and
$x+x'=r,$ \eqref{FDO.26} becomes
\begin{equation*}
\Id=r^{\ell+1}(1-s)^{\ell+2}\delta(s)\delta (\tilde y-\tilde y')
\delta (z-z')\beta_b^*(\nu'_\Phi).
\label{FDO.27}
\end{equation*}
Because $(1-s)\delta (s)=\delta (s)$, the factor $(1-s)^{\ell+2}$ may
be dropped.

Next, consider the lift from $X^2_b$ to $X^2_\Phi$. Of course,
the support of $\Id$ is contained in the $\Phi$-diagonal $D_\Phi$. 
In terms of the coordinates $r,$ $S=s/r,$ $Y=(\tilde y-\tilde y')/r$ 
valid near $D_\Phi $ (along with $\tilde y',z,z')$ we have
\begin{equation}
\Id=\delta (S)\delta (Y)\delta (z-z')\beta_\Phi^*(\nu '_\Phi).
\label{FDO.28}
\end{equation}
The density factor is simply a smooth, non vanishing,
section of the lift from the right factor of the density bundle $x^{-\ell-2}
dx\,dy\,dz;$ that is of the $\Phi$ density bundle.

For any embedded p-submanifold $M$ in a manifold with corners
$X$, the smooth $\delta $-functions on $M$ are the elements of a space
\begin{equation*}
\mathcal{D}^0(M)=\CI(X)\cdot\mu
\label{FDO.29}
\end{equation*}
where $\mu $ is any non-vanishing $\delta$-function with smooth
coefficients, as in \eqref{FDO.28}. The delta functions of order at most $k$
are obtained by differentiation
\begin{equation*}
\mathcal{D}^k(M)=\Diff k(X)\cdot\mu .
\label{FDO.30}\end{equation*}
In fact it is only necessary to differentiate across $M.$ Thus if
$\mathcal V$ is any Lie algebra of smooth vector fields which is transversal 
to $M,$ in the sense that for any section of $T_MX\big/TM=NM$ there is an 
element of $\mathcal V$ which projects to it along $M,$ then
\begin{equation*}
\mathcal{D}^k(M)=\sum\limits_{j\le k}\mathcal V^j\cdot \mu.
\label{FDO.31}
\end{equation*}

By Lemma~\ref{TLL}, $\mathcal{V}_\Phi $ lifts to $X^2_\Phi $ to
give such a Lie algebra transversal to $D_\Phi.$ Thus
\begin{equation}
\mathcal{D}^k(D_\Phi)\otimes \nu '_\Phi=\DiffF k\cdot\Id.
\label{FDO.32}
\end{equation}
The choice of initial smooth density on $M$ is of course
irrelevant; moreover, since it is in the right factor, the
differentiations on the right in \eqref{FDO.32} do not affect it. 
The space on the right here is, by definition, the 
space of Schwartz kernels of $\DiffF k(X),$ and so we conclude

\begin{prop}\label{CSK} The Schwartz' kernels of the elements of
$\DiffF k(X)$ lift to $X^2_\Phi $ to be precisely
the space $\mathcal{D}^k(D_\Phi)\cdot \nu '_\Phi $ of all smooth
$\delta$-functions on $D_{\Phi}$ up to order $k$ with a right
$\Phi$-density factor. 
\end{prop}

This result can be extended directly to operators on other
bundles. Thus if $E$ and $F$ are vector bundles over $X$ then
$\DiffF k(X;E,F)$ consists of all the differential operators from
sections of $E$ to sections of $F$ which are given by matrices with
elements in $\DiffF k(X)$ in local trivializations. It follows that
\begin{equation*}
\DiffF k(X;E,F)=\CI(X;F)\cdot\DiffF k(X)\cdot\CI(X;E^*).
\label{FDO.33}
\end{equation*}
The space on the right here is the finite linear span of (ordered)
products of elements from each of the three component spaces, hence
is simply the tensor product over $\CI(X)$. 
If $\Hom (E,F)$ is the bundle over $X^2$ with fibre $\hom(E_p,F_{p'})
= E_p^* \otimes F_{p'}$ at $(p,p')$ then it is also true that
\begin{multline}
\DiffF k(X;E,F)=\mathcal{D}^k(D_\Phi )\cdot
\CI(X^2_\Phi;\beta ^*_\Phi \Hom(E,F)\otimes\FDens')\\
=\mathcal{D}^k(D_\Phi )\cdot
\beta ^*_\Phi \CI(X^2_\Phi;\Hom(E,F)\otimes\FDens'),
\label{FDO.34}
\end{multline}
where $\FDens'$ is the lift of the $\Phi$-density bundle from the 
right factor. 

This gives the following normalization.

\begin{corollary}\label{CCSK} The lifts to $X^2_\Phi$ of Schwartz kernels 
of elements of $\DiffF k(X;E,F)$ coincides with the space 
\begin{equation}
\mathcal{D}^k(D_\Phi)\cdot\CI\left(X^2;\Hom(E;F)\otimes\FDens'\right).
\label{FDO.37}
\end{equation}
\end{corollary}

On a manifold without boundary this identification of the kernels of 
differential operators becomes 
\begin{equation*}
\Diff k(X';E.F)\equiv\mathcal{D}^k(D)\cdot
\CI\left(X^2;\Hom(E;F)\otimes\Dens'\right),
\label{FDO.38}
\end{equation*}
where again $\Dens'$ the density bundle lifted from the right factor.
To obtain the space of pseudodifferential operators in the boundaryless 
case, one replaces $\mathcal D^k(D),$ which is the space of polynomials
in all smooth vector fields, by $I^k(X^2,D),$ the space of conormal 
distributions, which may be thought of as symbolic functions of
these vector fields. Clearly 
\begin{equation*}
\mathcal D^k(D)\subset I^k(X^2, D), \label{FDO.39}
\end{equation*}
and in fact, $\mathcal{D}^k(D)$ may be characterized as that subspace 
of conormal distributions, whose elements have supports contained in $D.$ 

Although conormal distributions are initially defined with respect
to submanifolds of the interior, which do not intersect the boundary,
we may define conormal distributions with respect to any interior 
p-submanifold, simply by requiring that they extend across the 
boundary as conormal distributions for some (hence any) extension of the
submanifold. Thus $I^k(X^2_\Phi,D_\Phi ;G)$ is defined for any vector
bundle $G$ over $X^2_\Phi,$ and its elements are smooth outside $D_\Phi.$ 
Letting $\equiv$ denote equality in Taylor series we define the 
{\em microlocalization} of $\DiffF k(X)$ (or $\mathcal{V}_\Phi (x))$,
to be the following space of $\Phi$-pseudodifferential operators:

\begin{defin}\label{PPD} For any $m\in \RR$ the space of
$\Phi$-pseudodifferential operators of order $m$ (in the small calculus) is
\begin{multline}
\PsiF m(X;E,F)=\big\{K\in I^m(X^2_\Phi,D_\Phi;\beta ^*_\Phi 
\left(\Hom(E;F))\otimes\FDens'\right);\\
K\equiv0\text{ at } \pa X^2_\Phi \backslash\ff(X^2_\Phi )\big\},
\label{FDO.101}
\end{multline}
where $\ff(X^2_\Phi )$ is the front face produced by the blow up
\eqref{FDO.13}. 
There is some ambiguity in the definition of $I^m,$ depending on whether 
symbols of type $1,0$ or the smaller space of $1$-step polyhomogeneous
(\ie classical) symbols are used. When absolutely necessary, we shall denote
the polyhomogeneous space by 
\begin{equation*}
\PsipF m(X;E,F)\subset\PsiF m(X;E,F).
\label{FDO.40}
\end{equation*}
Generally the statements we make are valid with either interpretation of $I^m.$
\end{defin}

\section{Action of $\Phi$-pseudodifferential operators \label{S.AO}}

Combining Corollary~\ref{CCSK} and Definition~\ref{PPD} we have
\begin{equation}
\Dkf (X;E,F) \subset \PsiF k(X;E,F)\  \forall \ k\in \NN
\label{AO.1}\end{equation}
as spaces of kernels. We wish to interpret these Schwartz
kernels as operators so that \eqref{AO.1} still holds.

For simplicity take $E=F=\CC.$ Then the spaces in \eqref{FDO.37} and
\eqref{FDO.101} can be rewritten as
\begin{gather*}
\mathcal{D}^k(D_\Phi)\cdot\dCI_{\ff}(X^2;\beta^*_\Phi\FDens'))\\
I^m(X^2_\Phi,D_\Phi)\cdot\dCI_{\ff}(X^2;\beta^*_\Phi\FDens'),
\label{AO.2}
\end{gather*}
respectively, 
where $\dCI_{\ff}(X^2;\beta^*_\Phi\FDens')$ is the space of sections
vanishing to infinite order at all boundary faces except $\ff(X^2_\Phi).$

Consider the lift of a non-vanishing density from $X$ to the left factor of
$X^2$ and then to $X^2_\Phi.$ Using the diffeomorphism which exchanges
factors, the computation leading to \eqref{FDO.37} shows that the tensor
product identification $\rho_L^*\Dens\otimes\rho_R^*\FDens\equiv\Dens$
extends from the interior of $X^2,$ and so of $X^2_\Phi,$ to give an 
isomorphism of spaces of sections 
\begin{equation}
\dCI_{\ff}\left(X^2_\Phi;
\beta_\Phi^*(\rho_L^*\Dens\otimes\rho_R^*\FDens)\right)
\equiv\dCI_{\ff}(X^2_\Phi;\Dens).
\label{AO.101}\end{equation}
That is, the singular Jacobian factors arising all occur at faces other
than $\ff(X^2_\Phi)$ and have finite order singularities, which are absorbed
by the infinite order vanishing at these faces.

Fixing any $0<v\in\CI(X;\Dens)$, the action of $P\in\DiffF k(X)$ on 
$u\in\CI(X)$ can then be written 
\begin{equation*}
Pu\cdot v=(\pi_L)_*(P\cdot \pi_L^*v\cdot \pi_R^*u)
\label{AO.102}\end{equation*}
where \eqref{AO.101} is used to identify the product on the right as a
density on $X^2_\Phi$ and $\pi_L=\rho _L\circ\beta_{\Phi}$.
Generalizing from this, we see that in order to define the action of 
$\Phi$-pseudodifferential operators it suffices to establish 
the following result about push-forward: 

\begin{lemma}\label{PTC} Push-forward to the left factor defines a
continuous linear map 
\begin{equation*}
(\pi_L)_*:I^m(X^2_\Phi, D_\Phi)\cdot \dCI_{\ff}(X^2_\Phi;\Dens)
\longrightarrow\CI(X;\Dens).
\end{equation*}
and hence using \eqref{AO.101} 
\begin{multline}
(\pi_L)_*:I^m(X^2_\Phi, D_\Phi)\cdot \dCI_{\ff}(X^2_\Phi;
\beta_\Phi^*\pi^*_R\FDens)\longrightarrow\CI(X),\\
(\pi_L)_*(P\cdot\pi_L^*v)=(\pi_L)_*(P)\cdot v.
\label{AO.10}
\end{multline}
\end{lemma}

\begin{proof} This is a special case of the push-forward theorems discussed in
\cite{Melrose41}. To apply these theorems, we need to know that $\pi_L$ and 
$\pi_R$ are b-fibrations, and this is established below. The singularities of 
the kernel of $D_\Phi$ are integrated out since $\pi_L$ is a smooth map 
which is transversal to
the lifted diagonal. This transversality follows from Lemma~\ref{TLL} which
shows that the lift of $\mathcal{V}_\Phi(X)$ from the right factor is
transversal to $D_\Phi$ and spans the null space of the differential of
$\pi_L.$ From the general properties of conormal distributions
\begin{equation*}
(\pi_L)_*:I^m(X^2_\Phi, D_\Phi)\cdot\CIc(U;\Omega)\longrightarrow \CI(X)
\label{AO.11}
\end{equation*}
if $U$ is a small neighbourhood of the diagonal. For such a
neighbourhood $\CIc(U)\subset\dCI_{\ff}(X^2_{\Phi})$ so it suffices to
consider the case $m=-\infty,$ \ie to show that
\begin{equation}
(\pi_L)_*:\dCI_{\ff}(X^2_\Phi;\Dens)\longrightarrow\CI(X;\Dens).
\label{AO.12}
\end{equation}

Notice that \eqref{AO.12} is not quite trivial since it is {\em not} the
case that $(\pi_L)_*$ maps $\CI(X^2_\Phi;\Dens)$ into $\CI(X;\Dens).$
As discussed in \cite{Melrose41} a result such as \eqref{AO.12} follows from
two facts
\begin{equation*}
\begin{gathered}
\pi_L\text{ is a b-fibration and if } H\in M_1(X) \text{ then }\\
\pi^{-1}_L(H)\cap\ff\text{ is a boundary hypersurface of }X^2_\Phi.
\end{gathered}
\label{AO.13}
\end{equation*}
The second condition here just means that $f\in\dCI_{\ff}(X^2_\Phi)$
vanishes to infinite order on all of $\pi^{-1}_L(H)$ except for the one
boundary hypersurface, which is the front face corresponding to $H$ in
\eqref{FDO.13}.

Thus it is only necessary to show that $\pi_L=\rho_L\circ\beta _{\Phi}$
{\em is} a b-fibration (for the definition of this and other terms
here  we refer to \cite{Melrose41}).  Both $\rho_L$ and $\beta_\Phi$
are surjective b-maps and b-submersions, and so $\pi_L$ is also a b-submersion.
It remains only to see that no boundary hypersurface of $X^2_\Phi$ is mapped into a
boundary face of codimension two or more in $X$, but since $X$ does not have any such
faces this is automatically the case.
\end{proof}

Tensoring with the general coefficient bundle we deduce the elementary
mapping properties of $\Phi$-pseudodifferential operators.

\begin{prop}\label{DAO} Using the identification \eqref{AO.101}, each 
element $A\in \PsiF m(X;E,F)$ defines a continuous linear operator
\begin{equation}
A:\CI(X;E)\longrightarrow\CI(X;F)
\label{AO.14}
\end{equation}
which restricts to $A:\dCI(X;E)\longrightarrow\dCI(X;F)$ and extends by
continuity in the distributional topologies to
\begin{equation}
\begin{gathered}
A:\CmI(X;E)\longrightarrow \CmI(X;F) \Mand\\
A:\dCmI(X;E)\longrightarrow \dCmI(X;F).
\end{gathered}
\label{AO.15}
\end{equation}
These actions are consistent with the inclusion \eqref{AO.1}.
\end{prop}

\begin{proof} The discussion above proves \eqref{AO.14}, since any element can
be decomposed as a finite sum of products
\begin{equation*}
\Psi ^m_\Phi(X;E.F)\equiv\CI(X;F)\cdot \Psi ^m_\Phi (X)\cdot \CI(X;E^*).
\label{AO.16}
\end{equation*}
That $Au\in\dCI(X;F)$ if $u\in \dCI(X;E)$ follows from the observation that
if $\phi\in \CI(X)$ and $u\in \CI(X;E)$ then
\begin{equation*}
P(\phi u)=(\pi^*_R\phi\cdot P) u
\label{AO.17}\end{equation*}
where, directly from Definition~\ref{PPD}, $\pi^*_R\phi\cdot P\in
\Psi^m_\Phi(X;E,F)$ if $P\in \Psi ^m_\Phi(X;E,F).$ Now, if $u\in \dCI(X;E)$ it
can be written as a finite sum, $u=\sum _j\phi_j u_j, u_j\in \CI(X;E)$
and $\phi_j\in \dCI(X).$

If $\phi\in\dCI(X)$ then
\begin{equation*}
\pi^*_R\phi\cdot\dCI_{\ff}(X^2_\Phi )\subset \dCI(X^2_\Phi );
\end{equation*}
the extra vanishing at $\ff(X^2_\Phi )$ comes from the first factor.
From this it follows that
\begin{equation*}
\pi^*_L(x^{-k})\cdot \pi^*_R\phi P\in \Psi ^m_\Phi (X;E,F) \ \forall \ k \in \NN,
\quad \forall \ \phi \in \dCI(X),
\label{AO.18}
\end{equation*}
hence $x^{-k}Pu\in \CI(X;F),$ \ie $Pu\in \dCI(X;F).$

Again from Definition~\ref{PPD}, the formal adjoint of $P\in\Psi^m_\Phi(X;E,F)$ 
with respect to smooth inner products on $E$ and $F$ and a density on $X$ is 
an element of $\Psi ^m_\Phi (X;F,E).$ Thus the mapping
properties \eqref{AO.15} follow by duality.
\end{proof}

The singular function $x/x'$ on $X^2,$ where $x\in \CI(X)$  is a boundary
defining function on the left factor and $x'$ is the same function on the
right factor, lifts to be $\CI$ up to the interior of the front face of
$X^2_b,$ and hence up to the front face of $X^2_\Phi.$ Since it has only a
finite order singularity at the other boundary hypersurfaces, $\beta^*_\Phi 
(x/x')$ is a multiplier on $\Psi ^m_\Phi (X;E,F).$ This means in
particular that
\begin{equation}
\begin{gathered}
P_\pa:\CI(\pa X;E)\longrightarrow \CI(\pa X;F)\\
P_\pa u=P\tilde u\big|_{\pa X}, \tilde u\in\CI(X;E) \text{ with } u=\tilde
u\big|_{\pa X}\end{gathered}
\label{AO.19}
\end{equation}
is well defined, regardless of the extension $\tilde u$ of $u$. This corresponds 
to the map obtained by restricting an element of $\mathcal{V}_\Phi (X)$ to 
the boundary. Below it is augmented appropriately to define the normal 
operator, which is the boundary symbol in this context. Before doing this,
however, we first discuss the ordinary symbol map for $\Phi$-pseudodifferential 
operators.

For conormal distributions the symbol map
\begin{equation}
\begin{gathered}
I^m(X;G)\stackrel{\sigma _m}{\longrightarrow}S^{[m]}(N^*G;\Omega
^{\ha}(N^*G)\otimes\pi^*(\Omega ^{\ha}X)),\\
M=m+\frac 14 \dim X-\frac 12 \dim G
\end{gathered}
\label{AO.20}
\end{equation}
was introduced by H\"ormander. It is normalized on half-densities. Here
\begin{equation*}
S^{[M]}(\Lambda )=S^M(\Lambda )\big/ S^{M-1}(\Lambda )
\label{AO.21}
\end{equation*}
for any conic manifold $\Lambda ,$ is the quotient. For the case $G=D_\Phi
\subset X^2_\Phi $ it has already been shown in Lemma~\ref{TLL} that 
$N^*D_\Phi \simeq\FT^*X.$ The (singular) symplectic form on
$\FT^*X$ trivializes the bundle $\FDens=x^{-\ell-2}\Omega $
so \eqref{AO.20} leads to the desired map
\begin{equation}
\sigma _{\Phi ,m}:\PsiF m(X;E,F)\longrightarrow S^{[m]}(\FT^*X;\pi^*\hom(E,F)).
\label{AO.22}
\end{equation}
This generalizes the symbol map for differential operators obtained by
taking the leading part of \eqref{FCA.6} as a polynomial on $\FT^*X.$ It
gives a short exact sequence
\begin{equation*}
0\longrightarrow \PsiF{m-1}(X;E,F)\longhookrightarrow\PsiF m(X;E,F)
\stackrel{\sigma_{\Phi,m}}{\longrightarrow}S^{[m]}(\FT^*X;\pi^*\hom(E,F)).
\label{AO.23}
\end{equation*}
For the polyhomogeneous spaces, $\PsipF m(X;E,F)$ the symbol
becomes a homogeneous section of $\hom(E,F)$ lifted to $\FT^*X\backslash
0.$ Letting $Z=\FS^*X$ be the boundary `at infinity' 
of the radial compactification $\cFTs X$ this allows the symbol
map to be written
\begin{equation*}
\sigma_{p\Phi ,m}:\PsipF m(X;E,F)\longrightarrow
\CI(Z;(N^*Z)^m\otimes\pi^*\hom(E,F)),\ Z=\FS^*X.
\label{AO.24}
\end{equation*}

Next let us note how the action of $\Phi$-pseudodifferential operators
can be written locally.

\begin{prop}\label{AO.P} If $\chi\in\CI(X)$ has support in a coordinate
patch, $U,$ based at a boundary point $p\in\pa X$ with coordinates
$x,\tilde y,z$ as in Lemma~\ref{LCD} then the localized action of $P\in
\PsiF m(X)$ on $u\in\CIc(U)$ takes the form
\begin{equation}
\chi Pu=\int P_\chi (x,\tilde y,z,S,Y,z-z')\tilde v(x(1+xS),\tilde y-xY,z')
dS\,dY\,dz'\label{AO.25}
\end{equation}
where $\tilde v(x,y,z)$ is the coordinate representation of $u$ and the
kernel $P_\chi $ is the restriction to $U\times\RR^N$ of a
distribution on $\RR^n\times \RR^{n-k}\times\RR^k$ which has compact
support in the first and third variables, is conormal to
$\{S=0,Y=0\}\times\{z=z'\}$ (which  is the origin in the second two factors)
and is rapidly decreasing with all derivatives as
$|(S,Y)|\longrightarrow\infty.$ 
\end{prop}

\begin{proof} The kernel of the localized operator can be taken to be $\chi
P\chi .$ Any part of the kernel, on $X^2_\Phi ,$ away from $D_\Phi $ and
$\ff(X^2_\Phi )$ is smooth as a function on $X^2$ and vanishes rapidly at
the boundary. Localizing on $X^2_\Phi$ this gives a smooth section of
$\Hom(E,F)\otimes \pi^*_R\omega $ over  $X^2$ vanishing rapidly at both
boundaries; such a term can be written in the form \eqref{AO.25} with
$P_\chi $ both \ci\ and rapidly decreasing in $S$ and $Y.$

Thus we can suppose that the kernel has support in a small neighbourhood of
$D_\Phi \cup \ff(X^2_\Phi ).$ The part in the interior has a conormal
singularity at the diagonal and, since $x,x'\ne 0,$ can again easily be
written in the form \eqref{AO.25}. Thus we can suppose that the kernel has
support in a small neighbourhood of $\ff(X^2_\Phi ).$ Suppose initially that
its support only meets the interior of $\ff(X^2_\Phi ).$ In this region
$x,\tilde y,z,S=\frac{x-x'}{x^2}, Y=\frac{\tilde y-\tilde y}{x}$ and $z'$
gives coordinates on $X^2_\Phi .$ Thus \eqref{AO.25} results by introducing
$x'=x(1+Sx),y=\tilde y'-xY.$ The kernel has compact support and only a
conormal singularity at $S=0$, $Y=0,$ $z=z'$ so \eqref{AO.25} results.

The final term then is a smooth contribution to $P$ supported near
$\ff(X^2_\phi )$ and vanishing to infinite order at the other boundary faces
nearby. Although the coordinates $x,\tilde y, z,S,Y,z-z'$ are {\em not}
valid up to these adjacent boundaries a smooth function vanishing in Taylor
series in this sense just corresponds to a Schwartz function in the
variables $S,Y,$ \ie rapidly decreasing with all derivatives as
$|(S,Y)|\longrightarrow \infty.$ This proves the local representation
\eqref{AO.25}. 
\end{proof}

This proposition does not quite give a complete local description of the
action of $P.$ However, if $p,$ $ p'\in\pa X$ lie in the same fibre of
$\phi$ then $\phi(p)=\phi(p')$ lie in some coordinates patch in
$Y$. Thus one can take `consistent' coordinates near $p$ and $p'$ given by
$x,\tilde y,z$ and $x,\tilde y,\tilde z$ respectively. The same argument as
in the proposition gives a representation
\begin{equation}
\chi Pu=\int\tilde P_\chi (x,y,z,S,Y,\tilde z')\tilde v(x(1+xS),y-xY,\tilde
z')dS\ dY\ d\tilde z'
\label{AO.26}
\end{equation}
where $\chi \in \CI_c(X)$ has support in the coordinate patch near $p$ and
$v$ has support in the coordinate patch near $p'.$ The localized kernel
$\tilde P_\chi $ is smooth in all variables, compactly supported in $x,y,z$
and $\tilde z',$ and is rapidly decreasing with all derivatives as
$|(S,Y)|\longrightarrow \infty.$

Other pieces of the kernel correspond either to points $p,p'$ in different
fibres over the boundary or where either, or both, of the pair lie in the
interior. In these regions the localization of the kernel is a smooth section
of $X^2,$ except for a conormal singularity at the diagonal, and with
rapid vanishing at any boundary.

The front fact of $X^2_{\Phi}$ is a bundle over $\paX$ with fibre
$\phi^{-1}(\bar y)^2\times \cFN_{\bar y}T$ over $\bar y.$ The singular
variables $Y=(\tilde y-\tilde y')/x$ and $S=(x-x')/x^2$ introduced above 
give linear coordinates in $\FN_{\bar y}Y,$ depending on the choice of
admissible coordinates. Under a change of such coordinates $Y$ and $S$
transform linearly at $x=0,$ as a bundle transform on $\FN Y,$ and as Taylor
series at $x=0$ vary polynomially: 
\begin{equation*}
(Y,S)\longrightarrow A(\tilde
y')\cdot(Y,S)+\sum\limits_{j\ge1}x^jP_j(z,z',\tilde y',Y,S)
\label{AO.103}
\end{equation*}
where the $P_j$ smooth and are polynomials (without constant terms) in the
variables $Y,S.$ 

\section{Normal operator\label{S.NO}}

Using the representations \eqref{AO.25} and \eqref{AO.26}, we see that
the restriction map in \eqref{AO.19} is locally represented by
\begin{equation}
\begin{gathered}
P_\pa u=\int P_\pa(y,z,z-z')u(y,z')dz'\\
P_\pa(y,z,z-z')=\iint P_\chi (0,y,z,S,Y,z-z')dS\,dY.
\end{gathered}
\label{NO.1}
\end{equation}
This shows

\begin{lemma}\label{RRO} The map $P\longmapsto P_\pa$ in \eqref{AO.19} gives
a surjective map
\begin{equation}
\PsiF m(X;E,F)\longtwoheadrightarrow\Psif m(\pa X;E,F)
\label{NO.2}
\end{equation}
where $\Psif m(\pa X;E,F)$ is the space of pseudodifferential
operators acting on the fibres of $\phi:\pa X\longrightarrow Y$ and
depending smoothly on the base point.
\end{lemma}

It is important to note that the null space of \eqref{NO.2} consists of 
those elements for which the integral in \eqref{NO.1} vanishes for all 
$y\in Y$ (and $z, z' \in F$). This is closely related to the question
of determining which $\Phi$-pseudodifferential operators are
compact as operators on $L^2(X;E,F)$. For example, as will be seen below, 
the most obvious class of residual operators, the elements of 
$\PsiF{-\infty}(X;E,F)$, are not all compact. The operators $P$ for which
$P_\pa = 0$ are also not, in general, compact. In fact, such an operator
is compact only when the whole of the restriction of its kernel to 
$\ff(X^2_\Phi)$ vanishes, not just its fibre average as in \eqref{NO.1}. 

We examine this issue by means of `oscillatory testing'. To do this, 
fix a point $p\in\pa X,$ and suppose $f\in \CI(Y)$ is real-valued
and has $df(\phi(p))\ne 0.$ Choose $\tilde f\in \CI(X),$ also
real-valued, with $\tilde f\restrictedto\pa X=\phi^*f.$ Finally,
take $\chi\in\CI(X)$ such that 
\begin{equation*}
\begin{gathered}
\chi \equiv 1 \ \text{ near } \ \phi^{-1}(\phi(p))\\
d{\tilde f}\ne 0 \ \text{ on } \ \phi(\supp\chi \cap \pa X)
\end{gathered}
\label{NO.3}
\end{equation*}
and consider the `oscillatory test section'
\begin{equation*}
u_f=e^{i\tilde f/x}\chi u, \ u\in \CI(X;E).
\label{NO.4}
\end{equation*}

\begin{lemma}\label{OTN} For an `oscillatory test section' of this form,
and for any operator $P\in \PsiF m(X;E,F),$
\begin{equation}
P(e^{i\tilde f/x}\chi u)=e^{i\tilde f/x}\widetilde{P}u
\label{NO.5}
\end{equation}
with $\widetilde{P}\in\PsiF m(X;E,F).$
\end{lemma}

\begin{proof} The kernel of $\widetilde{P}$ is $e^{-i\tilde f/x}Pe^{i
\tilde f'/x'}\chi'$, using the obvious notation for variables 
and functions lifted from the left and right, respectively. It will 
suffice to show that 
\begin{equation}
\begin{gathered}
\text{the lift of } e^{-i\tilde f/x} \chi e^{i\tilde f'/x'}\chi' \text{ is
\ci\ on the union of } (X^2)^{\circ} \\ 
\text{and } (\ff(X^2_\Phi))^{\circ} 
\text{ and multiplication by it preserves } \dCI_{\ff}(X^2_\Phi).\end{gathered}
\label{NO.6}
\end{equation}
Recall that the space in the final statement here consists of
the smooth functions on $X^2_\Phi$ which vanish to infinite
order at all boundary faces except $\ff(X^2_\Phi)$.
The main point is to demonstrate the smoothness up to the interior of 
$\ff(X^2_\Phi ).$ First set $x'=(1+s)x,$ corresponding to the blow 
up \eqref{FDO.4}, so that  
\begin{equation*}
\frac{\tilde f}{x}-\frac{\tilde f'}{x'}=\frac{f(\tilde y)}{x}+g(x,\tilde
y,z)-\frac{f(\tilde y')}{x'}-\frac{g(x',\tilde y',z')}{1+s}.
\end{equation*}
Clearly we may restrict attention to the singular part
\begin{equation*}
\frac{f(\tilde y)}{x}-\frac{f(\tilde y')}{x'}
=\frac{f(\tilde y)}{x}-\frac{f(\tilde y-xY)}{(1+xS)x}.
\label{NO.7}
\end{equation*}
Using a Taylor expansion and the fact that $d_yf\ne0$ on $\supp\chi,$ 
we see that this is \ci\ up to $x=0$ as a function of $(S,Y)\in\RR^{1+k}.$ 
This proves the first part of \eqref{NO.6}; it also shows that this function 
has singularities only of finite order at all boundaries of $X^2_\Phi$ 
besides $\ff(X^2_\Phi ).$  It is therefore a multiplier on $\dCI_{\ff}(X^2_\Phi ).$ 
In fact, this shows that 
\begin{equation*}
e^{-i\tilde f/x}\chi e^{i\tilde f'/x'}\chi'\text{ is a multiplier on }
\Psi ^m_\Phi (X;E,F),
\label{NO.8}
\end{equation*}
and this proves the lemma.
\end{proof}

The operator $\widetilde{P}$ in \eqref{NO.5} depends not only on $f$, 
but also its extension $\tilde f$ and the cut off function $\chi.$ 
However, the restriction $\widetilde P_\pa$ depends only on $f$ and 
$\chi.$ This follows from \eqref{NO.1}, for if we assume that $\chi $ is 
supported in a coordinate patch, then 
\begin{multline}
\widetilde{P}_\pa(y,z,z-z')=\iint \exp(if'(y)\cdot Y+if(y)S)\\
P_\chi (0,y,z,S,Y,z-z')dS \, dY.
\label{NO.9}
\end{multline}
As noted in Lemma~\ref{RRO}, the $y$ variable enters here only as a parameter. 
A similar formula may be obtained from \eqref{AO.26}, and so we conclude 
that if $\chi (y)=1$ on the fibre $\phi^{-1}(y)$ then for that value of 
$y,$ $\widetilde P_\pa(y,z,z-z')\in\Psi ^m(\phi^{-1}(y);E,F)$ depends 
only on $f(y)$, $df(y)$ and $P.$

\begin{lemma}\label{LNO} If we fix a point $(y,\eta)\in T^*Y$ and a 
constant $\tau\in\RR$, then the indicial operator
\begin{equation*}
\widehat{P}(y,\tau,\eta)\in\Psi ^m(\phi^{-1}(y);E,F)
\label{NO.10}
\end{equation*}
is well-defined as the restriction to that fibre of $\widetilde P_\pa$,
where $\widetilde P$ is defined by Lemma~\ref{OTN} with $f$ chosen so that
$f(y)=\tau$ and $d_yf(y)=\eta.$  If $\widehat P(y,\tau,\eta )=0$ for every
$y,\tau ,\eta$, then $P\in x\Psi^m_\Phi (X;E,F).$
\end{lemma}

\begin{proof} Only the last statement needs to be checked. From
\eqref{NO.9} it follows that if $\widehat P(y,\tau ,\eta )\equiv 0$ then
the Fourier transform of the kernel on each fibre of $\ff(X^2_\phi)$ over
$Y$ vanishes, hence $P\restrictedto \ff(X^2_\Phi )\equiv 0$ and this is
equivalent to the existence of $Q\in \Psi ^m_\phi (X;E,F)$ such that $P=xQ$
(or equivalently, the existence of some $Q'$ with $P=Q'x').$
\end{proof}

As is clear from \eqref{NO.9}, \eqref{AO.25} and \eqref{AO.26} the
information carried in the operators $\widehat P$ as we let $y\in Y$, 
$f(y)=\tau$ and $\eta =df(y)$ vary determines the restriction of the kernel 
of $P$ to $\ff(X^2_\Phi).$ We shall reorganize these individual operators 
into the family of {\em normal operators.} Before we may do this, however,
we must first describe the algebra in which the normal family takes
values. 

For any compact manifold without boundary, $M$, and real vector space, $V,$
$M\times V$ is a $\CI$ manifold so the spaces $\Psi ^m(M\times V)$ of
pseudodifferential operators on $M\times V$ are well defined. These do not
compose since the growth of the kernels is unrestricted at infinity in $V.$
We consider the subspace
\begin{equation*}
\Psi ^m_{\sus(V)}(M)\subset\Psi ^m(M\times V)\label{NO.11}
\end{equation*}
consisting of the translation-invariant elements with $V$-convolution kernels
vanishing rapidly, with all derivatives at infinity. Thus
If $A\in\Psi ^m(M\times V)$ then $A\in \Psi ^m_{\sus(V)}(M)$ if
\begin{equation}
\begin{gathered}
AT^*_vu=T_v^* Au\ \forall \ u\in \CI_c(M\times V), \ v\in V\\
\text{ and } A: \CmIc(M\times V)\longrightarrow\CmIc(M\times V)+
\mathcal{S}(M\times V).
\end{gathered}
\label{NO.12}
\end{equation}
Here $T_v(m,w)=(m,w-v)$ is translation by $v$ and $\mathcal{S}(M\times V)$ is
the Schwartz space. The translation-invariance means that the kernel is of
the form
\begin{equation*}
A(m,m',v-v')\in\CmI(M^2\times V^2;\Omega _R).\label{NO.13}
\end{equation*}
Then, with some abuse of notation in which $A$ also stands for the
$V$-convolution kernel, the second condition in \eqref{NO.12} means that
\begin{equation}
A\in\CmI_c(M^2\times V; \Omega _R(M\times V))+\mathcal{S}(M^2\times V;\Omega
_R(M\times V))
\label{NO.14}
\end{equation}
where $\Omega _R(M\times V)=\pi^*_R\Omega,$ $\pi_R: \ M^2\times
V\longrightarrow M\times V$ being projection onto the right factor of $M.$

For the case $V=\RR$ this is the `suspended algebra'
considered in \cite{Melrose-Nistor2}. From \eqref{NO.14} and the general
properties of pseudodifferential operators it follows that $\Psi
^*_{\sus(V)}(M)$ is an order-filtered algebra of operators
\begin{equation*}
A: \ \mathcal{S}(M\times V)\longrightarrow \mathcal{S}(M\times V).
\label{NO.15}
\end{equation*}

The notation here, $\Psi ^*_{\sus(V)}(M),$ is to indicate that the algebra
can be thought of as the `$V$-suspended algebra of pseudodifferential
operators on $M.$' In this sense the primary object is $M.$ To have the
corresponding algebra of operators acting on a vector bundle, the vector
bundle should be defined over $M$ and pulled back to $M\times V.$ Thus if
$E$ is a bundle over $M$ then
\begin{equation*}
\Psi ^*_{\sus(V)}(M;E)=\Psi ^*_{\sus(V)}(M)\otimes_{\CI(M^2)}\CI(M^2;\Hom(E))
\label{NO.16}
\end{equation*}
defines the algebra of operators
\begin{equation*}
A: \ \mathcal{S}(M\times V;E)\longrightarrow \mathcal{S}(M\times V;E).\label{NO.17}
\end{equation*}

Directly from the definition, $\Psi ^m_{\sus(V)}(M;E)$ is invariant under
arbitrary diffeomorphism of $M$ and linear transformations of $V$, as well
as bundle transformations of $E$ over $M.$ This allows us to define the
more general object we need.

\begin{defin}\label{TIA} Let $\phi: X'\longrightarrow Y$ be a fibration
of compact manifolds, $E\longrightarrow X'$ a vector bundle and
$V\longrightarrow Y$ a real vector bundle. Then the algebra of
$V$-suspended fibre pseudodifferential operators on $X'$, $\Psi
^m_{\sus(V)-\phi} (X';E)$ is the space of operators
\begin{equation*}
A: \ \mathcal{S}(X'\times_YV;E)\longrightarrow
\mathcal{S}(X'\times_YV;E)\label{NO.18} 
\end{equation*}
which are local in $Y$ and for any open set $O\subset Y$ over which $\phi$ and
$V$ are trivial reduce to a smoothly $O$-parametrized element of $\Psi
^m_{\sus(V_y)}(\phi^{-1}(y);E).$
\end{defin}

Thus an element  $A\in \Psi ^m_{\sus(V)_\phi}(X';E)$ has Schwartz
kernel of the form
\begin{gather*}
A(y,z,z',v)\in \CmI_c(X'\times_Y X'\times_YV;\Hom(E)\otimes \Omega _R)\\
+\mathcal{S}(X'\times_Y X'\times_Y V;\Hom(E)\otimes\Omega _R)
\end{gather*}
where $A$ is conormal with respect to the submanifold
\begin{equation*}
D_\Phi \times\{0\}=\{(y,z,z,0)\}\subset X'\times_YX\times_YV
\label{NO.19}
\end{equation*}
which is the fibre diagonal. The action of $A$ is given explicitly by
\begin{equation*}
Au(y,z,v)=\int A(y,z,z', v-v')u(y,z',v')dz' dv'
\label{NO.20}
\end{equation*}
since $\Omega _R$ is the lift from the right factor of the fibre density
bundle of $X'\times_YV$ as a fibration over $Y.$

Since the kernel is essentially a density on the fibres of $V$ when all
the variables are held fixed its Fourier transform is well defined and is a
smooth function of the dual variables
\begin{equation*}
\hat A(y,z,z',w^*)=\int e^{-iw^*\cdot w} A(y,z,z',w)dV.
\label{NO.21}
\end{equation*}
For each $w^*\in V^*_y$ it is a pseudodifferential operator on the fibre
$\phi^{-1}(y).$ This corresponds to the indicial operator in Lemma~\ref{LNO}.
In fact

\begin{prop}\label{TNH} For a $\Phi$-structure on a compact manifold with
boundary $X,$ the indicial operators of Lemma~\ref{LNO} combine to give the
Fourier transform of an element of $\Psi ^m_{\sus(V)-\phi}(\pa X;E,F)$
where $V=\FN Y$ is the null bundle, on $Y$, of the restriction ${}^\Phi
T_{\pa X} X\longrightarrow T\pa X$ and the resulting map, defining the
normal operator, gives a short exact sequence
\begin{equation}
0\longrightarrow x\Psi ^m_\Phi (X;E,F)\longrightarrow \Psi ^m_\Phi (X;E,F)
\stackrel{N_\Phi }{\longrightarrow} \Psi ^m_{\sus(V)-\phi}(\pa
V;E,F)\longrightarrow 0.
\label{NO.22}
\end{equation}
\end{prop}

\begin{proof} From \eqref{NO.9} we know that $\widehat P$ is the Fourier
transform of the restriction of the kernel of $P$ to the front fact,
$\ff(X^2_\Phi).$  Thus, at the level of kernels, the map $N_\Phi $ is just
restriction to $\ff(X^2_\Phi ).$ This shows that the null space of $N_\Phi $
acting  on $\Psi ^m_\Phi (X;E,F)$ is precisely $x\Psi ^m_\Phi (X;E.F)$ and
directly from Definition~\ref{TIA}, $N_\Phi $ is surjective as is
\eqref{NO.22}.
\end{proof}

When we consider composition below it will be apparent that \eqref{NO.22}
is multiplicative.

\section{Composition\label{S.C}}

It is relatively straightforward, if tedious, to check that the space
$\Psi^*_\Phi(X;E)$ is an algebra by using the local representations
\eqref{AO.25} and \eqref{AO.26}. Instead we use a more conceptual approach
that has the virtue of applying in rather general circumstances
\cite{Melrose8} and in the present circumstances to more general operators
(\ie `larger calculi' with non-trivial boundary behaviour).

Thus our approach is to define a `triple $\Phi $ product' $X^3_\Phi $ with
maps back to the double product $X^2_\Phi $ defined in \eqref{FDO.15}.
The definition of $X^3_\Phi $ from $X^3$ proceeds by a chain of five blow ups.
These are carried out independently at each of the boundary faces of $X,$
so for simplicity we generally assume that $\pa X=H$ is connected. We shall
use a notation for the boundary faces of $X^3$ similar to that used above
for $X^2.$ Namely if $H\in M(X)$ then set
\begin{equation*}
L(H)=H\times X^2,\ M(H)=X\times H\times X,\ R(H)=X^2\times H.
\label{C.1}
\end{equation*}
Thus in general
\begin{equation*}
M_1(X^3)=\bigcup\limits_{H\in M_1(X)}\{L(H), \ M(H), \ R(H)\}.
\label{C.2}
\end{equation*}
For the codimension two boundary faces we are only interested in those
meeting the diagonal; we use the notation
\begin{equation*}
S(H)=H\times H\times X, \ C(H)=H\times X\times H, \ F(H)=X\times H\times H.
\label{C.3}
\end{equation*}
Here `$S=$  second', `$C=$ composite' and `$F=$ first' arise from the
relationship to the composition of operators. The only codimension three
boundary faces meeting the diagonal are
\begin{equation*}
T(H)=H^3\subset M_3(X^3),\label{C.4}
\end{equation*}
the `triple' faces. In general we drop the reference to $H.$

The two stage blow up leading to $X^3_b$ resolves the intersection of $T,$
$S,C$ and $F$:
\begin{equation}
X^3_b=[X^3; T; S; C; F].\label{C.5}
\end{equation}
Although there are in principle four blow ups here, after the blow up of
$T$ the lifts of $S,C$ and $F$ are disjoint p-submanifolds so can be
blown up in any order.

The remaining stages in the definition of $X^3_\Phi$ involve the blow up of
various $\Phi$-diagonal submanifolds. To see how these arise, consider the
product
\begin{equation*}
X\times X^2_b=[X^3; F].\label{C.6}
\end{equation*}
The submanifold $\Phi \subset X^2_b$ defined in \eqref{FDO.10} therefore
defines a submanifold we denote
\begin{equation}
\Phi_F=X\times \Phi \subset X\times X^2_b.
\label{C.7}
\end{equation}
Now $T\subset F$ so, by the commutativity of blow-up in this setting (see
\cite{Melrose8}, \cite{Epstein-Melrose-Mendoza1}), the order of
blow ups can be exchanged to obtain a natural isomorphism
\begin{equation*}
[X^3;F;T]\simeq [X^3;T;F].\label{C.8}
\end{equation*}
The product structure in \eqref{C.7} and the fact that $T$ lifts
to $[X^3;F]$ to be
\begin{equation*}
T'=H\times B_b(H)\subset [X^3;F]=X\times X^2_b,\label{C.9}
\end{equation*}
shows that $\Phi _F$ has a common product decomposition with $T'.$ The
inverse image of $\Phi _F$ in $[X^3;F;T]$ is therefore the union of two
p-submanifolds which we denote
\begin{equation*}
\begin{cases}
\widetilde\Phi _F=\beta ^*_{T'}(\Phi _F)=\clos((\beta _{T'})^{-1}(\Phi
_F\backslash T')) \text{ and}\\
\widetilde\Phi _{FT}=(\beta _{T'})^{-1}(\Phi_F\cap T').
\end{cases}
\label{C.10}
\end{equation*}
Neither of these p-submanifolds meets the lifts of $S$ or $C$ to
$[X^3;F,T]$ so they equally well define submanifolds
\begin{equation*}
\widetilde{\Phi}_F,\widetilde{\Phi }_{FT}\subset X^3_b\equiv [X^3;F;T;S;C].
\label{C.11}
\end{equation*}
Of course from the basic symmetry of the set up we have similar submanifolds
\begin{equation*}
\widetilde{\Phi }_S, \widetilde{\Phi }_{ST}, \widetilde{\Phi }_C,
\widetilde{\Phi }_{CT}\subset X^3_b.
\label{C.12}
\end{equation*}
Notice that $\widetilde{\Phi }_O\subset O_b\subset X^3_b, $ $O=F,S,C$ when
$O_b$ denotes the front face produced by the blow up of $O$ in defining
$X^3_b.$ On the other hand $\widetilde{\Phi}_{OT}\subset T_b$ for $O=F,S,C.$

\begin{lemma}\label{TFD} The intersection of any pair of $\widetilde{\Phi
}_{ST},$ $\widetilde{\Phi }_{FT}$ and $\widetilde{\Phi }_{CT}$ is the
submanifold
\begin{equation*}
\widetilde{\Phi }_{T}=\widetilde{\Phi }_{ST}\cap \widetilde\Phi _{FT}\cap
\widetilde{\Phi }_{CT}\label{C.13}
\end{equation*}
which is contained in the interior of $T_b.$
\end{lemma}

\begin{proof} Let us examine these definitions more closely. Since we only
need to consider the operations near each boundary, $X$ can be replaced by
$[0,1)_x\times H,$ so $X^2_b$ is given by \eqref{FDO.18} and in this
representation
\begin{equation*}
\Phi =\{0\}\times\{0\}\times D_\Phi\label{C.14}
\end{equation*}
where $D_\Phi\subset H\times H$ is the fibre diagonal. Now, $X^3\simeq
[0,1)^3\times H^3.$ So near the new faces
\begin{equation*}
\begin{gathered}{}
[X^3;T]\simeq [0,1)\times G\times H^3\\
X^3_b\simeq [0,1)\times G_b\times H^3.
\end{gathered}
\label{C.15}
\end{equation*}
Here, $G\subset \RR^2$ is an equilateral triangle with centre the origin
and $G_b$ is obtained by blowing up each corner of it. Thus $G_b$ can be
embedded in $\RR^2$ as a regular hexagon with centre the origin. The sides
of this hexagon are alternately the front faces and original boundaries,
\ie $C_b, R_b, F_b, M_b, S_b, L_b.$


The lifts of the $\Phi$  diagonals are easily identified, thus
\begin{equation}
\widetilde \Phi _F\simeq [0,1)\times \{p_F\}\times H\times D_\Phi
\label{C.16}
\end{equation}
where $p_F\in G_b$ is the midpoint of the side corresponding to $F_b$ and
$D_\Phi\subset H^2$ is the $\phi$-fibre diagonal. Similarly


\begin{equation*}
\widetilde\Phi _{FT}=\{0\}\times \ell_F \times H\times D_\Phi
\label{C.116}
\end{equation*}
where $\ell_F\subset G_b$ is the line through $p_F,$ the origin and the
midpoint of the side representing $L_b.$

This proves the lemma with
\begin{equation}
\widetilde\Phi _T=\{0\}\times \{0\}\times T_\Phi, \ T_{\Phi}\subset
H^3 \text{ the triple } \Phi\text{-diagonal}.
\label{C.17}
\end{equation}

Now we complete the definition of the triple $\Phi$-space by three more
(levels of) blow up
\begin{equation}
X^3_\Phi =[X^3_b;\widetilde\Phi _T;\widetilde\Phi _{FT};\widetilde\Phi_{ST};
\widetilde\Phi _{CT}; \widetilde\Phi _F;\widetilde\Phi _S;\widetilde\Phi _C].
\label{C.18}
\end{equation}
From \eqref{C.16} and \eqref{C.17} it follows that $\widetilde\Phi _{FT},
\widetilde\Phi _{ST};\widetilde\Phi _{CT}$ lift to be disjoint after the
blow up of $\widetilde\Phi _T$ so the orders of these three blow ups, and
the last three, are immaterial. However, the order between the last three
blow ups and the preceding three is important and cannot be arbitrarily
rearranged, since for instance $\widetilde\Phi _{FT}$ and $\widetilde\Phi
_F$ intersect but not transversally, nor is one contained in the
other. This space is mainly useful for the maps defined on it.
\end{proof}

\begin{prop}\label{DPM} For $O=F,S,C$ there is a b-fibration
$\pi^3_{\Phi ,0}:X^3_\Phi \longrightarrow X^2_\Phi$ fixed by the demand
that it give a commutative diagramme with the corresponding projection
\begin{equation}
\xymatrix{
X^3 \ar[r]^{\psi^3}\ar[d]^{\pi^3_{\Phi,0}}
&X^3_b\ar[r]^{\beta^3_b}\ar[d]^{\pi^3_{b,0}}
&X^3\ar[d]^{\pi^3_0}\\
X^2_\Phi\ar[r]_{\psi ^2} &X^2_b \ar[r]_{\beta^2_b}& X^2.}
\label{C.19}
\end{equation}
\end{prop}

\begin{proof} To define these maps we start with the corresponding maps for
the b-calculus; the middle maps in \eqref{C.19}. These can be constructed
using the commutability of blow ups for $O\supset T,$ we shall take $O=F$
for the sake of definiteness. Then
\begin{equation}
X^3_b=[X^3;T;F;S,C]=[X^3;F;T;S,C]=\bigl[[X^3; F]; T;S,C\bigr].
\label{C.20}
\end{equation}
Now $[X^3;F]=X\times X^2_b$ so there is a commutative diagram with vertical
projections.
\begin{equation*}
\xymatrix{
X\times  X^2_b\ar[r]\ar[d]^{\gamma_F}&X^2\ar[d]^{\pi^3_F}\\
X^2_b\ar[r]^{\beta ^2_b}&X^2.}
\label{C.21}
\end{equation*}
Then $\pi^3_{b,F}=\gamma _F\circ \tilde \beta ,$ $\tilde\beta:
\bigl[[X^3;F];T;S;C\bigr]\longrightarrow [X^3;F]$
being the blow down map. Thus $\pi^3_{b,F}$ is defined and is automatically
a b-map. We need to show that it is a b-submersion and finally a
b-fibration. Certainly it is surjective.

A b-fibration, $f,$ remains a b-submersion when composed with the blow
down map for blow up of some p-submanifold, $M,$ if, for each point p of the
submanifold the induced map
\begin{equation}
f:M \longrightarrow\Fa(f(p))
\label{C.22}
\end{equation}
is a b-submersion. Here $\Fa(q)$ is the smallest boundary face of the
range space containing $f(p).$ For any boundary face, $M,$ this condition
is automatically satisfied. This `blown up' b-fibration is again a
b-fibration, rather than just a b-submersion, if $f(M)$ is a boundary
hypersurface of the range space, which is to say it is not contained in a
boundary face of codimension $2.$ Since this is immediately clear for the
blow ups is the definition of $\tilde\beta ,$ and hence $\pi^3_{b,F},$ the
latter map is a b-fibration. 

Now that we have fixed the central vertical maps in \eqref{C.19} we proceed to
the definition of the $\pi^3_{\Phi ,0},$ again taking $O=F$ for
definiteness sake. In \eqref{C.18} the submanifolds $\widetilde{\Phi
}_{ST},$ $\widetilde{\Phi}_{CT}$ and $\widetilde\Phi _F$ are disjoint, so the
order can be changed to
\begin{equation*}
X^3_\Phi =[X^3_b; \widetilde\Phi_T; \widetilde\Phi_{FT}; \widetilde\Phi_F;
\widetilde\Phi'],\
\widetilde\Phi '=\widetilde\Phi _{ST};\widetilde\Phi _{CT};\widetilde\Phi
_S;\widetilde\Phi _C. 
\label{C.23}
\end{equation*}
Similarly $\widetilde\Phi_T\subset \widetilde\Phi_{FT}$ and
$\widetilde\Phi_T$ is disjoint from $\widetilde\Phi_F$ so
\begin{equation}
X^3_\Phi =[X^3_b; \widetilde\Phi _{FT};\widetilde\Phi _F;\widetilde\Phi
''],\ \widetilde\Phi ''=\widetilde\Phi _T;\widetilde\Phi '.
\label{C.24}
\end{equation}
Consider again the definition, \eqref{C.5}, of $X^3_b,$ reorganized as in
\eqref{C.20}. The submanifold $S$ and $C$ are disjoint from $\widetilde\Phi
_{FT}$ and $\widetilde\Phi _F$ so \eqref{C.24} can be written
\begin{equation}
X^3_\Phi =[X\times X^2_b;T;\widetilde\Phi _{FT};\widetilde\Phi _F;R],
\ R=S;C;\widetilde\Phi ''.\label{C.25}
\end{equation}
In $X\times X^2_b$ the submanifold $X\times \Phi $ lifts to $\widetilde\Phi
_F$ under the blow up of $T$ and $\widetilde\Phi _{FT}$ is the lift, in
fact preimage, of $\phi (X\times\widetilde\Phi )\cap T$ under blow up of
$T.$ Thus \eqref{C.25} can be commuted to
\begin{equation*}
X^3_\Phi =\bigl[X\times X^2_b; (T\cap (X\times\Phi ));T;\widetilde\Phi
_F;R\bigr]. \label{C.26}
\end{equation*}
The second and third blow up are disjoint so in fact
\begin{equation}
X^3_\Phi =\bigl[X\times X^2_b; (X\times\Phi ); \ T\cap (X\times\Phi
);T;R\bigr]
\label{C.27}
\end{equation}
The final rearrangement here is of two cleanly intersecting submanifolds
with are blown up with there intersection, this can be accomplished by
blowing up either of them first, then the intersection, then the other,
with the same final result.

The first blow up in \eqref{C.27} is the definition of $X^2_\Phi $ so
\begin{equation}
X^3_\Phi =[X\times X^2_\Phi ;T\cap (X\times\Phi ); T;R]
\label{C.28}
\end{equation}
allows the blown up projection in \eqref{C.19} to be defined by
\begin{equation*}
\pi^3_{\Phi ,F}=\tilde\gamma _F\cdot \tilde\psi , \tilde\gamma _F:X\times
X^2_\Phi \longrightarrow X^2_\Phi \label{C.29}
\end{equation*}
being the projection, with $\tilde\psi $ the collective blow up of $R$ in
\eqref{C.28}.

To show that $\pi^3_{\Phi ,F}$, and hence by symmetry each of the
$\pi^3_{\Phi ,0},$ is a b-fibration it is only necessary to check the two
conditions involving \eqref{C.22} for each of the blow ups in $\tilde\psi.$
In fact, using the product structure in \eqref{C.16}, etc., this is
straightforward so the details are omitted. Suffice it to say that the
fibration can be eliminated directly and the case $\phi=\Id$ is then
simpler to analyze.
\end{proof}

We further augment Proposition~\ref{DPM} by considering the relationship
between these maps and the lifted diagonals.

\begin{lemma}\label{TFDS} The lifted diagonals, defined as the closures in
$X^3_\Phi $ of the diagonal $D\subset X^{\circ}\times X^{\circ}$ in each of
the three possible positions, are p-submanifolds $D_{\Phi ,F},$ $D_{\Phi
,S}$, $D_{\Phi ,C}$ as is the lifted triple diagonal $D_{\Phi ,T}$. Each of
the maps $\pi^3_{\Phi,O}$ is transversal to $D_{\Phi ,O'}$, for $O'\neq O$
and maps
\begin{equation*}
D_{\Phi ,T}=D_{\Phi ,O_1}\cap D_\Phi ,O_2, O_1\neq O_2
\label{C.30}
\end{equation*}
diffeomorphically onto $D_\Phi \in X^2_\Phi.$
\end{lemma}

\begin{proof} These results are immediate away from any boundaries. The
transversality of $\pi^3_{\Phi ,F}$, say, to $D_{\Phi ,S}$ follows by
lifting $\mathcal{V}_\Phi(X) $  from the left factor. This is in the null
space of 
the differential of $\pi^3_{\Phi ,F}$ and lifts to be transversal to
$D_{\Phi ,S},$ essentially by Lemma~\ref{TLL}. Thus $\pi^3_{\Phi ,F}$ maps
$D_{\Phi ,S}$ diffeomorphically onto $X^2_\Phi$ and hence embeds the
submanifold $D_{\Phi ,T}\subset D_{\Phi ,S}$ as $D_{\phi }\subset X^2_\Phi.$
\end{proof}

With these maps and transversality results available the composition formula
is now straightforward.

\begin{theorem}\label{MCT} For any vector bundles $E,$ $F,$ $G$ over a compact
manifold with boundary $X$, and fibred boundary structure $\Phi ,$
\begin{equation*}
\Psi ^m_\Phi (X;F,G)\circ\Psi ^{m'}_\Phi (X;E,F)\subset\Psi ^{m+m'}_\Phi
(X;E,G)
\label{C.31}
\end{equation*}
and both the symbol map \eqref{AO.22} and normal operators
\begin{equation*}
\Psi ^m_\Phi (X;E.F)\longrightarrow \Psi ^m_{\sus(V)-\phi}(\pa X;E,F),
\label{C.32}
\end{equation*}
$V={}^\Phi NY,$ of Proposition~\ref{TNH}, are multiplicative.
\end{theorem}

\begin{proof} The composition is well defined by Proposition~\ref{DAO}.
\end{proof}

\section{Mapping properties\label{S.MP}}

To deduce the $L^2$ boundedness of the operators of order zero we shall use
an argument due to H\"ormander \cite{Hormander6} which depends on the
existence, within the calculus, of an approximate square root of a positive
elliptic element.

\begin{prop}\label{SQT} If $B\in\PsiF 0(X)$ is formally self-adjoint,
for some smooth positive density on $X,$ then for $C>0$ sufficiently large 
\begin{equation*}
C+B=A^*A+R,
\label{ML.1}\end{equation*}
for some $A\in\PsiF 0(X)$ and $R\in x^{\infty}\PsiF {-\infty}(X).$
\end{prop}

\begin{proof} Since $B$ is formally self-adjoint with respect to the density,
$\nu,$ the indicial family $\widehat B(\tau,\eta)$ consists of operators
which are self-adjoint with respect to the boundary density, defined by
$\nu =dx\otimes\nu _0$ for an admissible defining function $x.$ Thus, for
$C>0$ sufficiently large 
\begin{equation*}
(C+\widehat B(\tau,\eta))^{\ha}\in\Psif0(\pa X)
\label{ML.2}\end{equation*}
and from the uniqueness of this positive square root it is the indicial
family of some $A_0\in\PsiF 0(X).$ Again for $C$ large enough $A_0$
can be chosen to have 
\begin{equation*}
\sigma_0(A_0)=(C+\sigma_0(B))^{\ha}
\label{ML.3}\end{equation*}
as well. Thus, replacing $A_0$ by $\ha(A_0+A_0^*)$ we find 
\begin{equation*}
C+B-A_0^2\in x\PsiF {-1}(X).
\label{ML.4}\end{equation*}

Proceeding by induction, as in the standard case, one can suppose that
$A_{(k-1)}\in\PsiF 0(X)$ has been constructed such that
$A_{(k-1)}^*=A_{(k-1)}$ and 
\begin{equation*}
C+B-A_{(k-1)}^2=R_k\in x^k\PsiF {-k}(X).
\label{ML.5}\end{equation*}
Adding an unknown $A_k\in x^k\PsiF {-k}(X)$ to $A_{(k-1)}$ gives 
\begin{equation*}
\begin{gathered}
C+B-(A_{(k-1)}+A_k)^2=R_k-A_{(k-1)}A_k-A_kA_{(k-1)}-A_k^2\\
\equiv R_k-A_{(k-1)}A_k-A_kA_{(k-1)}
\end{gathered}\label{ML.6}\end{equation*}
modulo $x^{k+1}\PsiF {-k-1}(X).$ Thus if $A_k=x^kG_k$ is chosen to satisfy
\begin{equation}
N(A_0)N(G_k)+N(G_k)N(A_0)=N(F_k),\ F_k=x^{-k}R_k
\label{ML.7}
\end{equation}
then $A_{(k)}=A_{(k-1)}+A_k$ satisfies the inductive hypothesis at the next
level. Notice that, at the level of the indicial families, \eqref{ML.7} is
indeed solvable, as the linearization of the definition of the square root
\begin{equation*}
(\widehat{A_0}(\tau,\eta)+\widehat{G_k}(\tau,\eta))^2=
\widehat{A_0}(\tau,\eta)^2+\widehat{F_k}(\tau,\eta),
\label{ML.8}\end{equation*}
$\widehat{A_0}(\tau,\eta)$ being a positive operator for all $\tau,\eta.$
Finally then $A$ can be taken as an asymptotic sum of the series defined by
the $A_k.$
\end{proof}

\begin{theorem}\label{ML.9} Each element $P\in\PsiF 0(X;E)$ defines
a bounded linear operator on $L^2(X;E),$ defined with respect to a positive
smooth density on $X.$
\end{theorem}

\begin{proof} Since $X$ is compact, boundedness on $L^2$ is a local property of
operators, so it suffices to consider the case $E=\CC$ by local
trivialization. Then applying Proposition ~\ref{SQT} with $B=-P^*P$ shows
that, for all $u\in\dCI(X),$
\begin{equation*}
\|Pu\|^2 = C\|u\|^2-\|Au\|^2+\langle Ru,u\rangle \le C\|u\|^2+
|\langle Ru,u\rangle|\le C'\|u\|^2,
\label{ML.10}\end{equation*}
where the fact that elements of $x^{\infty}\PsiF {-\infty}(X),$ being
smoothing operators, are $L^2$ bounded has been used.
\end{proof}

Just as the construction of an approximate square root proceeds as in the
boundaryless case, with some extra care needed to handle the normal
operator, so the existence of parametrices for `fully elliptic' operators is
straightforward.

\begin{prop}\label{ML.12} If $P\in\PsiF m(X;E,F)$ is {\em fully
elliptic} in the sense that its symbol is everywhere invertible and its
normal operator is invertible on each fibre of $\phi,$ then there exits
$Q\in\PsiF {-m}(X;F,E)$ satisfying  
\begin{equation*}
P\circ Q-\Id\in x^{\infty}\PsiF {-\infty}(X;F)\Mand Q\circ P-\Id\in
x^{\infty}\PsiF {-\infty}(X;E).
\label{ML.13}\end{equation*}
\end{prop}

\begin{proof} Using the symbol calculus, $Q_0\in\PsiF {-m}(X;F,E)$ can
be chosen to have
\begin{equation*}
\sigma_{-m}(Q_0)=(\sigma_m(P))^{-1},\ N(Q)=N(P)^{-1}.
\label{ML.14}\end{equation*}
This ensures that $P\circ Q_0=\Id-R_1,$ with $R_1\in x\PsiF {-1}(X;F).$
Proceeding inductively it can be supposed that $Q_j\in
x^j\PsiF {-m-j}(X;F,E)$ have been constructed so that 
\begin{equation*}
P\circ(\sum\limits_{j=0}^{k-1}Q_j)=\Id-R_kx^k,\ R_k\in
\PsiF {-k}(X;F).
\label{ML.15}\end{equation*}
Adding $Q_k=T_kx^k\in x^k\PsiF {-m-k}(X;F,E)$ where
$\sigma_{-m-k}(T_k)=(\sigma_m(P))^{-1}\sigma_{-k}(R_k)$ and
$N(T_k)=N(P)^{-1}N(T_k)$ gives the next inductive step. Then $Q$ can be
taken to be an asymptotic sum of the $Q_k.$
\end{proof}

As in the boundaryless case these basic results easily lead to
continuity, compactness and Fredholm properties on Sobolev spaces. For
positive real number $m,$ and any $l\in\RR$ set 
\begin{equation*}
\begin{gathered}
x^lH_{\Phi}^m(X;E)=
\left\{u\in x^lL^2(X;E);Pu\in L^2(X;E)\ \forall\ P\in\PsiF m(X;E)\right\}\\
{\begin{aligned}x^l&H_{\Phi}^{-m}(X;E)\\
&=\left\{u\in\CmI(X;E);u=\sum\limits_{i=1}^N P_iu_i,\ u_i\in x^lL^2(X;E),\
P_i\in\PsiF m(X;E)\right\}\end{aligned}}
\end{gathered}
\label{ML.11}\end{equation*}

\begin{lemma}\label{ML.18} For these $\Phi$-Sobolev spaces 
\begin{equation*}
x^lH_{\Phi}^m(X;E)\subset x^{l'}H_{\Phi}^{m'}(X;E)\Longleftrightarrow l\ge
l'\Mand m\ge m'
\label{ML.19}\end{equation*}
with the inclusion then continuous. The inclusion is compact if and only if
$l>l'$ and $m>m'$ and each $P\in\PsiF m(X;E,F)$ defines a continuous
linear map
\begin{equation}
P:x^lH_{\Phi}^{m'}(X;E)\longrightarrow x^lH_{\Phi}^{m'-m}(X;F)
\label{ML.21}
\end{equation}
for all $l$ and $m'.$
\end{lemma}

\begin{prop}\label{ML.17} Each fully elliptic element,
$P\in\PsiF m(X;E,F),$ is Fredholm as a map \eqref{ML.21} and
conversely this condition characterizes the fully elliptic elements. The null
space of such an operator is contained in $\dCI(X;E)$ and there is
a complement to the range in $\dCI(X;F).$
\end{prop}

\begin{prop}\label{ML.16} If $P\in\PsiF m(X;E,F)$ is fully elliptic
then $P^*P+1$ has a two-sided inverse in $\PsiF {-2m}(X;E).$
\end{prop}

\section{Wavefront set\label{S.WS}}

There is a natural notion of wavefront set associated to the calculus of
operators $\PsiF*(X;E).$ In fact in a certain sense there are two such
notions, one associated to regularity and the other associated to growth at
the boundary. In each case we first consider the corresponding notion of
microlocal support, or operator wavefront set, for the operators before
examining the wavefront set of distributions.

For an embedded submanifold $Y$ of a manifold $X$ the conormal distributions
introduced by H\"ormander, $I(X,Y),$ have wavefront set a closed conic
subset of the conormal bundle to $Y$ in $X.$ Let $SN^*Y$ be the boundary of
the compactification of this bundle, \ie the quotient of $N^*Y\setminus0$ by
the $\RR^+$-action. Then
\begin{equation*}
\WF(u)\subset SN^*Y,\ u\in I^*(X,Y)
\end{equation*}
can also be identified with the cone support of the symbol obtained by
transverse Fourier transformation of $u.$ This second definition
extends directly to the case of an interior p-submanifold of a manifold
with corners. In particular it applies to the lifted diagonal in
$X^2_{\Phi}.$ This allows us to define the `symbolic' part of the
$\Phi$-wavefront set by 
\begin{gather*}
\WFFs'(A)=\WF(A)\subset SN^*(\Diag_{\phi})=\FS^*X,\\
\WFFs'(A)=\emptyset\Longleftrightarrow A\in\PsiF{-\infty}(X;E). 
\end{gather*}
The elliptic subset, $\Ells^m(A)\subset\WFFs'(A)$ is the open subset of
$\FS^*X$ on which the symbol of order $m$ has an inverse of order $-m.$
Here we have used the identification of the conormal bundle to the lifted 
diagonal with $\FT^*X.$

Now, the discussion above of the composition of $\Phi$-pseudodifferential
operators shows that the diagonal singularity of the composite arises from
the same operation as in the interior case. In particular the standard
proof of the microlocality of composition shows that 
\begin{equation}
\WFFs'(A\circ B)\subset\WFFs'(A)\cap\WFFs'(B),\ A,B\in\PsiF*(X;E).
\label{WS.15}
\end{equation}
The construction of parametrices for elliptic operators can also be
microlocalized, so if $K\subset\Ells^m(A)$ is closed, for a given $A\in\PsiF
m(X;E),$ then there exists $B\in\PsiF{-m}(X;E)$ such that 
\begin{multline}
A\circ B=\Id - R_L,\ B\circ A=\Id - R_R,\ R_L,\ R_R\in\PsiF0(X;E)\Mand\\
K\cap\big(\WFFs'(R_L)\cup\WFFs'(R_R)\big)=\emptyset.
\label{WS.18}
\end{multline}

Combining these standard results extended to the $\Phi$-calculus leads to
an alternative characterization of the operator wavefront set

\begin{lemma}\label{WS.16} For any $A\in\PsiF*(X;E)$ 
\begin{equation}
\big(\WFFs'(A)\big)^{\complement}=\bigcup\big\{\Ells^0(B);B\in\PsiF0(X;E)\Mand
B\circ A\in\PsiF{-\infty}(X;E)\big\}.
\label{WS.17}
\end{equation}
\end{lemma}

\begin{proof} If $p\in\FS^*X$ is in the set on the right in \eqref{WS.17}
then there is some $B\in\PsiF0(X;E)$ which is elliptic at $p$ and such
that $B\circ A\in\PsiF{-\infty}(X;E).$ Using a microlocal parametrix as in
\eqref{WS.18} it follows that $p\notin\WFFs'(A).$ The converse inclusion
follows from the microlocality, \eqref{WS.15}.
\end{proof}

We next define the corresponding notion of support, $\WFFs(u),$ for any
distribution $u\in\CmI(X).$ Since operators of order $-\infty$ are ignored
here we work modulo the space 
\begin{equation*}
x^{-\infty}H^{\infty}_{\Phi}(X)=
\bigcup\limits_{k\in\ZZ}x^kH^{\infty}_{\Phi}(X),
\end{equation*}
Indeed, 
\begin{equation*}
A\in\PsiF{-\infty}(X)\Longrightarrow
A:\CmI(X)\longrightarrow x^{-\infty}H^{\infty}_{\Phi}(X).
\label{WS.38}\end{equation*}
Then we simply define 
\begin{gather*}
\WFFs(u)=\bigcap\left\{\Char_{\Phi}(A);
A\in\PsiF0(X),\ Au\in x^{-\infty}H^{\infty}_{\Phi}(X)\right\}\subset\FS^*X,\\
\Char_{\Phi}(A)=(\Ells^0(A))^{\complement}.
\end{gather*}

Thus by definition, $p\notin\WFFs(u)$ if there exists $A\in\PsiF0(X)$ which
is elliptic at $p$ and is such that $Au\in
x^{-\infty}H^{\infty}_{\Phi}(X).$ As with the standard wavefront set there
is an alternate characterization in terms of the essential support 
\begin{multline}
\left(\WFFs(u)\right)^{\complement}=
\bigcup\big\{U;U\subset\FS^*X\text{ is open s.t. }\\
A\in\PsiF0(X),\ \WFFs'(A)\subset U
\Longrightarrow Au\in x^{-\infty}H^{\infty}_{\Phi}(X)\big\}.
\label{WS.39}
\end{multline}
This follows by use of the calculus as in the boundaryless case. From
\eqref{WS.39}, or directly, the calculus is microlocal for this wavefront
set: 
\begin{equation*}
\WFFs(Au)\subset\WFFs'(A)\cap\WFFs(u),\ A\in\PsiF{\infty}(X),\ u\in\CmI(X).
\label{WS.40}\end{equation*}
Note also that 
\begin{equation*}
u\in\CmI(X)\Mand\WFFs(u)=\emptyset\Longrightarrow
u\in x^{-\infty}H^{\infty}_{\Phi}(X).
\label{WS.59}\end{equation*}

Together with this extension of the usual notion of wavefront set
we next consider related notions at the boundary. First consider the operator
wavefront set. This will be defined as a subset of the radial
compactification $\cFNs\pbX$ of the bundle $\FN^*\pbX.$ This `$\Phi$-conormal
bundle' to the fibres of the boundary is the space of parameters in the
normal operators; note that it is a bundle over $\pbX,$ the base of the
fibration, and that it is the dual of the bundle corresponding to the Lie
subalgebra in \eqref{FCA.100}. Its lift to $\paX,$
$\FN^*\paX=\phi^*(\FN^*\pbX),$ occurs as the quotient of the part,
$\FT^*_{\paX}X,$ of the dual of the structure bundle over the boundary by
the subbundle  
\begin{equation}
\fT^*\paX=\bigcup\limits_{p\in\paX}T^*\phi^{-1}(p)\subset\FT^*_{\paX}X,\
\fpi:\FT^*_{\paX}X\longrightarrow \FN^*\paX,
\label{WS.67}
\end{equation}
of the fibre cotangent bundles. The inclusion here is just given by pairing
with vector fields, which shows $\fT^*\paX$ to be the annihilator bundle
in $\FT^*_{\paX}X$ of the lift of $\FN\pbX.$

Now, let $\bar y\in\pbX$ be a point in the base of the fibration of the
boundary and consider a finite point $p\in\FN^*_{\bar y}Y.$ For any
admissible coordinates $x,$ $\tilde y$ near $\bar y,$
$p=d(\frac{\bar\lambda +\bar\eta\cdot(\tilde y-\bar y)}x)$ for some
$\bar\lambda,$ $\bar\eta.$

Then we \emph{define}
\begin{multline}
p\notin(\WFFb'(A)\cap\FN^*Y)\Longleftrightarrow
\fS^*_{\bar y}X\cap\WFFs'(A)=\emptyset\Mand\exists\ \psi \in\CI(X)\Mst\\
\exp{\left(-i\frac{\lambda+\eta\cdot(\tilde y-\bar y)}x\right)}\psi A
\exp{\left(i\frac{\lambda+\eta\cdot(\tilde y-\bar y)}x\right)}:
\CI(X)\longrightarrow\dCI(X)\\
\forall\ (\lambda ,\eta )
\text{ in some neighbourhood of }(\bar\lambda ,\bar\eta ),
\label{WS.32}
\end{multline}
where $\psi\in\CI_\phi (X),$ is of the form
$\psi=\phi^*\psi'$ on the boundary with $\psi'(\bar y)=1$ and $\psi'$ is
supported in the coordinate patch.

Thus in order that $p\notin\WFFb'(A)$ we first demand that $\WFFs'(A)$ not
meet $\fS^*_{\phi^{-1}(\bar y)}X.$ Note that the preimage of $p$ in
$\FN^*_{\paX}X$ under projection to $\FN^*\paX$ and then $\FN^*\pbX$ meets
the sphere bundle at infinity $\FS^*{\paX}X$ exactly in $\fS^*_{\phi
^{-1}(\bar y)}X.$ Thus this is the condition that the part of $\WFFs'(A)$
`lying above' $p$ should be trivial. The second part of \eqref{WS.67}
implies in particular that the normal operator of $A$ should be trivial
near $p.$ In fact, in terms of the local representation \eqref{AO.25} and
\eqref{AO.26}, it means that the Fourier transform in $S$ and $Y$ of the
local kernel should vanish in a fixed neighbourhood of the point $(\lambda
,\eta)$ and $\tilde y'=\bar y$ as a function $z,z'$ and in the sense of
Taylor series in $x.$ The uniformity of the neighbourhood in $x$ is
important. It follows from the remarks after \eqref{AO.26}, in particular
the polynomial dependence of the coordinate transformation, that this
condition, of vanishing, is independent of coordinates.

Thus this notion is independent of the choice of coordinates in
\eqref{WS.32}. It is clearly multiplicative in the usual sense that
\begin{equation*}
\WFFb'(AB)\subset\WFFb'(A)\cap\WFFb'(B)\Mif A,B\in\PsiF*(X).
\end{equation*}
There is an importance difference between operators of finite order and
those of order $-\infty$ as regards $\WFFb'.$ Of course, for the latter the
condition on $\WFFs'(A)$ in \eqref{WS.32} is vacuous and then given a point
$p\in\FN^*\pbX$ with a neighbourhood $U$ we can always decompose
\begin{multline}
\PsiF{-\infty}(X)\ni B=B'+B'',\ B',B''\in\PsiF{-\infty}(X),\\
p\notin\WFFb'(B''),\ \WFFb'(B')\subset U,\ p\in\FN^*\pbX.
\label{WS.46}\end{multline}
Such a decomposition is not in general possible for operators of
finite order, since for instance the ellipticity of the symbol of $B$
would imply that the indicial operator never vanishes.

To an infinite point $p\in\FSN^*_{\bar y}\pbX$ there corresponds a `preimage'
$\Gamma(p)\subset\FS^*_{\bar y}X,$ consisting of the intersection 
\begin{equation}
\Gamma(p)=\clos\left(\fpi^{-1}(p')\right)\cap\FS^*X\Min \cFTs X.
\label{WS.54}
\end{equation}
Here $p'\subset\FN^*Y$ is the open half line corresponding to the point $p$
on the sphere at infinity and $\ftpi$ is the composite of $\fpi$ in
\eqref{WS.67} and the projection from $\paX$ to $\pbX.$ Thus $\Gamma(p)$
is a closed half-sphere bundle of fibre dimension $\dim F +1$ over $\phi
^{-1}(\pi(p)).$ We define the condition
\begin{multline}
p\notin\WFFb'(A)\Mfor p\in\FSN^*Y\Longleftrightarrow
\Gamma(p)\cap\WFFs'(A)=\emptyset\Mand\\
(\gamma\cap\FN^*\pbX)\cap \WFFb'(A)=\emptyset\text{ for some open
}\gamma\subset\cFNs\pbX\Mwith p\in\gamma.
\label{WS.55}\end{multline}
Note that \eqref{WS.32} shows that the analogue of $\Gamma(p)$ in case
$p\in\FN^*{\bar y}Y$ is finite is $\fS^*_{\phi^{-1}(\bar y)}X.$ If
$p\FSN^*{\bar y}Y$ then $\Gamma(p)\supset\fS^*_{\phi^{-1}(\bar y)}X.$ 

The restriction of the conjugated operator  
\begin{equation*}
\exp{\left(-i\frac{\lambda+\eta\cdot(\tilde y-\bar y)}x\right)}\psi A
\exp{\left(i\frac{\lambda+\eta\cdot(\tilde y-\bar y)}x\right)}
\end{equation*}
in \eqref{WS.32} to the boundary fibre above $\pi(p)$ is the indicial
operator, $N(A,p),$ at $p.$ We define ellipticity for operators of
order $m$ in this boundary sense by 
\begin{multline}
\Ellb^m(A)=\{p\in\FN^*\pbX;N(A,p)^{-1}\text{ exists in
}\Psi^{-m}(\phi^{-1}(\pi(p))\}\\
\cup\{p\in\FSN^*\pbX;\Gamma(p)\subset\Ells^m(A)\}\subset\cFNs\pbX.
\label{WS.52}
\end{multline}
Then certainly $\Ellb^m(A)\subset\WFFb'(A).$

\begin{lemma}\label{WS.63} For any $A\in\PsiF m(X)$ the set $\Ellb^m(A)$ is
open in $\cFNs\pbX.$
\end{lemma}

\begin{proof} Certainly if $p\in\Ellb^m(A)\cap\FN^*\pbX$ then $\Ellb^m(A)$
contains a neighbourhood of $p,$ since the invertibility of the normal
operator is an open condition. So consider $p\in\FSN^*\pbX$ `at infinity'
and suppose $p\in\Ellb^m(A).$ Since $\Gamma(p')$ is compact and depends
continuously on $p'\in\FSN^*\pbX$ it follows that
$\Gamma(p')\subset\Ells(A)$ for $p'$ in a neighbourhood of $p.$ Thus it
remains to show that $N(A,q)^{-1}\in\Psi^{-m}(\phi^{-1}(\pi(q)))$ for
$q\in\gamma'\cap\FN^*\pbX$ for some neighbourhood $\gamma '$ of $p$ in
$\cFNs\pbX.$ Using the calculus, we may construct an operator
$G\in\PsiF{-m}(X)$ such that $G\circ A=\Id-E$ where
$\Gamma(p)\cap\WFFs'(E)=\emptyset.$ Shrinking $\gamma'$ as necessary, it
follows that $N (E,q)$ is in $\Psi^{-\infty}(\phi^{-1}(\pi(q)))$ for
$q\in\gamma '\cap\FN^*\pbX$ and is rapidly vanishing as $q\to\FSN^*\pbX$
in $\gamma '.$ Thus $N(A,q)^{-1}\in\Psi^{-m}(\phi^{-1}(\pi(q)))$ exists for all
$q$ in the intersection of $\FN^*\pa X$ and some neighbourhood of $p$ in
$\cFNs\pa X.$ Thus $\Ellb^m(A)$ is open. 
\end{proof}

The construction in the proof of this lemma can be slightly extended to
yield:

\begin{lemma}\label{WS.66} If $p\in\cFN^*\pbX$ and $A\in\PsiF{m}(X)$
then $p\in\Ellb^m(A)$ if and only if there exists $G\in\PsiF{-m}(X)$ such
that $p\notin\WFFb'(\Id-A\circ G),$ $p\notin\WFFb'(\Id-G\circ A).$ 
\end{lemma}

Notice that in demanding that $A$ be elliptic at a finite point
$p\in\FN^*\pbX$ we are requiring that $A$ be symbolically elliptic on the
whole set $\fS_{\pi(p)}^*\paX\subset\FS^*_{\pi(p)}\paX,$ which is the
sphere of the subspace in \eqref{WS.67} above the point $\pi(p),$ since
$N(A,p)$ is to be invertible as a pseudodifferential operator of order $m$
on the boundary fibre. Correspondingly if $p\in\FN^*_{\bar y}Y$ then the
parametrix $G$ may be chosen to have $\WFFs'(G)$ concentrated near
$\fS_{\phi^{-1}(\bar y)}^*\paX\subset\FN^*X$ whereas $\WFFb'(G)$ can only be
concentrated near the fibre $\cFN^*_{\bar y}\pbX.$ If $p\in\FSN^*_{\bar
y}Y$ then $\WFFs'(G)$ may be concentrated near $\Gamma(p)$ and again
$\WFFb'(G)$ may be concentrated near $\cFN^*_{\bar y}\pbX.$

We now define
\begin{multline}
\WFFb(u)^\complement=\big\{p\in\cFNs\pbX;\exists\ A\in\PsiF0(X),\
p\in\Ellb^0(A),\ Au=w+\sum\limits_jB_jv_j,\\ 
w\in\dCI(X),\ v_j\in\CmI(X),\ B_j\in\PsiF{-\infty}(X),\
p\notin\WFFb'(B_j)\big\}\subset\cFNs\pbX.
\label{WS.41}
\end{multline}
Taking 
\begin{equation*}
\Charb^m(A)=\cFNs\pbX\setminus\Ellb^m(A)
\label{WS.56}\end{equation*}
this can also be written 
\begin{multline}
\WFFb(u)=\bigcap\bigg\{\Charb^0(A)\cup\bigcup\limits_j\WFFb'(B_j);
A\in\PsiF0(X),\ B_j\in\PsiF{-\infty}(X)\Mwith\\
Au=w+\sum\limits_jB_jv_j,\Mforsome w\in\dCI(X),\ v_j\in\CmI(X)\bigg\}.
\label{WS.64}
\end{multline}
The extra finite sum of terms $B_jv_j$ is included in \eqref{WS.41}, and
\eqref{WS.64}, because
of the non-localizability of $\WFFb'$ for operators of finite order. Notice
that if $B\in\PsiF{-\infty}(X)$ has $\WFFb'(B)$ concentrated sufficiently
close to $p\notin\WFFb(u),$ so $\WFFb'(B_j)\cap\WFF'(B)=\emptyset$ for each
$j,$ then $B\sum_jB_jv_j\in\dCI(X)$ too.

Since we are demanding that $Au$ lie in the `residual space' at $p$
\begin{multline}
\mathcal{R}_p(X)=\big\{u\in\CmI(X);u=u_1+u_2,\ u_1\in\dCI(X),\\
u_2\in\{B\in\PsiF{-\infty}(X);\ p\notin\WFFb'(B)\}\cdot\CmI(X)\big\}\subset
x^{-\infty}H^{\infty}_{\Phi}(X),
\label{WS.58}\end{multline}
where the $\cdot$ means finite span; this is a considerably finer
notion than $\WFFs(u)$ over the boundary.

\begin{lemma}\label{WS.44} If $p\in\cFN^*\pbX,$ the condition
  $p\notin\WFFb(u)$ given by \eqref{WS.41} is equivalent to the existence
  of $u_1\in\dCI(X),$ $C\in\Psi^0(X)$ with $p\notin\WFFb'(C),$
  $u'_j\in\CmI(X)$ and $B_j\in\PsiF{-\infty}(X)$ for $j=1,\dots,J$ with
  $p\notin\WFFb'(B_j)$ such that
\begin{equation}
u=u_1+\sum\limits_jB_ju'_j+Cu.
\label{WS.45}
\end{equation}
\end{lemma}

\begin{proof} The form \eqref{WS.45} for $u$ follows by applying the
parametrix $G$ of Lemma~\ref{WS.66} to the defining relation in \eqref{WS.41}.

Conversely, if \eqref{WS.45} holds for $p\in\FN^*\pbX$ and $A\in\PsiF0(X)$
is elliptic at $p$ and has $\WFFs'(A)$ in a small neighbourhood of
$\fS_{\pi(p)}^*\paX,$ so that $\WFFs'(A)\cap\WFFs'(C)=\emptyset$ then 
\begin{equation*}
  Au=Au_1+\sum\limits_jAB_ju'_j+ACu\in\mathcal{R}_p(X),
\label{WS.A2}
\end{equation*}
since $\WFFs'(AC)=\emptyset.$ This gives \eqref{WS.41}. A similar argument
applies if $p\in\FSN^*Y.$
\end{proof}

As already noted, the subtlety with the definition of $\WFF(u)$ above
arises from the non-localizability of the normal operators. In the
particular case of the scattering calculus, considered in \cite{Melrose43}
and \cite{Melrose-Zworski1}, there is no such difficulty. It is
useful to relate the general case to this scattering case. 

\begin{lemma}\label{WS.60} If $\psi \in\CI(X)$ has support sufficiently
close to $\phi ^{-1}(\bar y)\subset\paX$ for some point $\bar y\in Y$ then
there is an open product neighbourhood of $\supp(\psi)$ of the form  
\begin{equation*}
[0,\epsilon)_x\times Y'\times F,\ Y'\subset Y,
\end{equation*}
consistent with the fibration of the boundary and then for any
$A\in\PsiF{-\infty}(X),$ $\psi A\psi$ is a smooth right density on $F\times F$
with values in the scattering calculus on $X'=[0,1]\times Y,$ that is
$\Psisc{-\infty}(X').$ Furthermore, this product decomposition allows
$\FN^*\pbX$ to be identified with $\cscT^*_{Y}X'$ and if $B\Psisc{0}(X')$
is supported sufficiently close to the boundary and has
$\WFsc'(B)\cap\cscT^*_{Y}X'\cap\WFFb'(C)=\emptyset,$ where
$C\in\PsiF{0}(X)$ is supported in the product neighbourhood then
$B\circ C\in\rho^{\infty}\PsiF{-\infty}(X').$
\end{lemma}

\begin{proof} The first part follows directly from the definitions of the
algebras in terms of their kernels on the blown up spaces since locally, in
$Y,$ the blow up defining the stretched product for the fibred cusp
calculus is just the blow up for the scattering calculus (\ie the case that
the fibres in the boundary are points) with the fibres $F\times F$ as
factors.

The composition statement in the second part follows directly from the local
normal forms \eqref{AO.25}, \eqref{AO.26}.
\end{proof}

Despite the complexity of its definition, we may now see that this notion of
wavefront set has many of the familiar properties.

\begin{prop}\label{WS.48} The set $\WFF(u)=\WFFs(u)\cup\WFFb(u)\subset \FS^*
  X\cup\cFNs Y$ is closed, is empty only for elements of $\dCI(X),$
  satisfies
$$
\WFFb(u_1+u_2)\subset\WFF(u_1)\cup\WFF(u_2)
$$ and
  is reduced by the application of pseudodifferential operators,
  $A\in\PsiF*X,$ in the sense that
\begin{equation*}
\WFF(Au)\subset\WFF'(A)\cap\WFF(u),\ \WFF'(A)=\WFFs'(A)\cup\WFFb'(A).
\label{WS.49}\end{equation*}
\end{prop}

\begin{proof} That $\WFF(u)$ is closed follows directly from the
openness of the elliptic sets. The microlocality of pseudodifferential
operators, \eqref{WS.49}, follows directly for the interior part of the
wavefront set and from \eqref{WS.45} for the boundary part. Thus, if
$B\in\PsiF*X$ and $p\notin\WFFb(u)$ then first applying $A$ to
\eqref{WS.45} and then applying $Q\in\PsiF0X$ which is elliptic at $p$ but
has small support (see the discussion following Lemma~\ref{WS.66}) gives 
\begin{equation*}
QAu=QAu_1+\sum\limits_jQAB_ju_j+QACu.
\end{equation*}
Here in the last term, $QAC\in\PsiF{-\infty}X$ if $\WFFs'(Q)$ is chosen
sufficiently small and $p\notin\WFFb'(QAC).$ Thus it can be absorbed as an
extra term in the sum and deduce that $p\notin\WFFb'(Au)$ by
\eqref{WS.64}. The other components of \eqref{WS.49} are simpler.

It remains to show that if $\WFF'(u)=\emptyset$ then $u\in\dCI(X).$ From
$\WFFs(u)=\emptyset$ it follows that $u\in
x^{-\infty}H^{\infty}_{\Phi}(X);$ in particular it is smooth in the
interior of $X.$ We may localize the support of $u$ to a small set near a
boundary point, using the microlocality just discussed; thus we may assume
that $u$ has small support, in which the fibration has a product
decomposition. Thus $u(x,\tilde y,z)$ is a smooth function of $z$ with
values in a fixed space $x^{-N}H^{-N}_{\text{sc}}(X'),$ $X'=[0,1)_x\times
Y$ as in Lemma~\ref{WS.60}. Applying the second half Lemma~\ref{WS.60}, it
follows that if $A\in\Psisc0{X'}$ has wavefront set concentrated near any
point $p\in\cscT^*_{Y}X'$ then, applying it to \eqref{WS.45} $Au(x,y,z)$ is
\ci\ in $z$ with values in $\dCI(X'),$ and hence in $\dCI(X'\times F).$
Applying this to a partition of unity in the scattering calculus it follows
that $u\in\dCI(X).$ 
\end{proof}

\begin{remark}\label{WS.51} The somewhat global (at least on the fibre)
condition in \eqref{WS.41}, coming in turn from \eqref{WS.52}, is
necessitated by the fact, mentioned above, that one cannot freely localize
the indicial family. Thus, if $A\in\PsiF0(X)$ has indicial family
invertible, in the calculus, at any one point $p\in\FN_yY$ its indicial
family cannot be zero at any other point in that fibre, that is, 
\begin{equation*}
p\in\FN^*_{\bar y}X,\ p\in\Ellb^0(A)\Longrightarrow \cFNs_{\bar
y}\subset\WFFb'(A).
\end{equation*}
\end{remark}

\section{Fibred cusp metrics\label{S.FCM}}

As an application of the discussion above we shall examine the
spectrum of the Laplacian for a metric of `exact $\Phi$-type'. By this we
mean any Riemann metric on the interior of $X,$ a manifold with a fibred
boundary as in \eqref{I.1}, which takes the form 
\begin{equation}
g=\frac{dx^2}{x^4}+\frac{h'}{x^2}+g'
\label{FCM.1}
\end{equation}
for some product decomposition near the boundary
$X\subset[0,\epsilon)_x\times\pa X$ with $g_Y$ a smooth symmetric
$2$-cotensor on $[0,\epsilon)\times Y$ which is positive definite when
restricted to $\{0\}\times Y$ (with restriction $h)$ and $g'$ is a smooth
symmetric $2$-cotensor on $X$ which is positive definite when restricted to
each fibre over the boundary. The fibration $\phi$ and the boundary
defining function $x$ in \eqref{FCM.1} together determine a $\Phi$
structure on $X.$ Moreover

\begin{prop}\label{FCM.2} The Laplacian of a metric \eqref{FCM.1} is a
$\Phi$-differential operator on functions or acting on sections of the
$\Phi$ exterior bundle.
\end{prop}

The metric $g$ is a positive definite metric on the bundle $\FT X,$ smooth
and non-degenerate up to the boundary. This allows $\fT^*\paX$ to be
identified with the orthocomplement of $\FN^*\paX$ in $\FT^*_{\pa X}X.$
Furthermore, the boundary defining function $x$ in \eqref{FCM.1} defines a
natural section $dx/x^2$ of $\FN^*\paX$ the orthocomplement of which can be
identified with the lift of $T^*Y,$ by identifying $\eta\cdot dy$ with
$\frac{\eta\cdot dy}x.$ For each $y\in\paX$ let $\Lap_y$ be the
Laplacian on the fibres $\phi^{-1}(y)$ fixed by the metric $g'.$ Let
$\lambda_j(y)$ be the eigenvalues of $\Lap_y$ arranged in increasing order,
repeated with multiplicity.

\begin{theorem} If $u\in\CmI(X;\Lambda^k)$ satisfies $\Lap u-\lambda
u\in\dCI(X),$ with $\lambda \in\CC$ then 
\begin{gather}
\lambda \notin[0,\infty)\Longrightarrow u\in\dCI(X),\label{FCM.3}\\
{\begin{split}\lambda \in[0,\infty)\Longrightarrow
\WFF(u)\subset\big\{q\in\FN_y\paX;\ \exists\ \lambda_j(y)\le\lambda\Mst\\ 
q=s\frac{dx}{x^2}+\frac{\eta}x\Mwith s^2+|\eta|_h^2=\lambda-\lambda_j(y)\big\}.
\end{split}}
\label{FCM.4}\end{gather}
\end{theorem}

\ifx\undefined\bysame
\newcommand{\bysame}{\leavevmode\hbox to3em{\hrulefill}\,}
\fi

\end{document}